# SEMIPARAMETRIC ESTIMATION OF FRACTIONAL COINTEGRATING SUBSPACES

By Willa W. Chen[1] and Clifford M. Hurvich

*Texas A&M University and New York University*

We consider a common-components model for multivariate fractional cointegration, in which the $s \geq 1$ components have different memory parameters. The cointegrating rank may exceed 1. We decompose the true cointegrating vectors into orthogonal fractional cointegrating subspaces such that vectors from distinct subspaces yield cointegrating errors with distinct memory parameters. We estimate each cointegrating subspace separately, using appropriate sets of eigenvectors of an averaged periodogram matrix of tapered, differenced observations, based on the first $m$ Fourier frequencies, with $m$ fixed. The angle between the true and estimated cointegrating subspaces is $o_p(1)$. We use the cointegrating residuals corresponding to an estimated cointegrating vector to obtain a consistent and asymptotically normal estimate of the memory parameter for the given cointegrating subspace, using a univariate Gaussian semiparametric estimator with a bandwidth that tends to $\infty$ more slowly than $n$. We use these estimates to test for fractional cointegration and to consistently identify the cointegrating subspaces.

**1. Introduction.** Fractional cointegration has been the subject of much recent attention; see, for example, the work of Robinson [16], Robinson and Marinucci [19], Robinson and Marinucci [18], Chen and Hurvich [3]. All of these papers assume either that the observed series is bivariate or that the cointegrating rank is 1. Arguably the most interesting case, from an econometric point of view, is the situation where the series is multivariate and has cointegrating rank which may exceed 1. This situation was covered by Robinson and Yajima [20], who considered methods for determining the cointegrating rank, and also by Chen and Hurvich [4], who focused on estimation of the space of cointegrating vectors.

Received November 2003; revised January 2006.
[1]Supported by NSF Grant DMS-03-06726.

Chen and Hurvich [4] studied the properties of eigenvectors of an averaged periodogram matrix of differenced, tapered observations, averaging over the first $m$ Fourier frequencies, where $m$ is held fixed as the sample size grows. They showed that the eigenvectors corresponding to the $r$ smallest eigenvalues (where $r$ is the cointegrating rank) lie close to the space of true cointegrating vectors with high probability. They also presented an empirical analysis of fractional cointegration in US interest rates for bonds of seven different maturities. They found evidence that the cointegrating rank was greater than one and, furthermore, that the memory parameter of the cointegrating errors may take on a variety of values that differ substantially if cointegrating vectors corresponding to substantially different eigenvalues are used. This last finding, while of apparent interest from an econometric point of view, could not be explained directly from the theoretical results presented in [4] since they did not attempt in their theory to separate the space of cointegrating vectors into subspaces yielding different memory parameters.

The goals of the present paper are to exhibit a model that allows us to highlight these subspaces, to show that the subspaces and their corresponding memory parameters can be estimated individually and to show how to use the residual-based Gaussian semiparametric estimates of the memory parameters to consistently identify the cointegrating subspaces and to test for fractional cointegration. By contrast, Chen and Hurvich [4] did not consider either testing for cointegration or estimation of the degree of cointegration.

We first present, in Section 2, a semiparametric common-components model in which the components have different memory parameters, while the entries of the observed multivariate series have just one common memory parameter. Next, we show that the space of cointegrating vectors can be decomposed into a direct sum of orthogonal *cointegrating subspaces* such that vectors from distinct subspaces yield cointegrating errors with distinct memory parameters.

In Section 5, we show that each of these cointegrating subspaces can be separately estimated using sets of eigenvectors of the averaged periodogram matrix. Since $m$ is held fixed, we are able to obtain a rate of convergence for the estimated cointegrating vectors that depends only on the difference between the memory parameters in the given and adjacent subspaces and is not hampered by the rate of increase of $m$, as in other related work (cf. [19], in the bivariate case).

To each true cointegrating subspace, there corresponds an estimated cointegrating subspace spanned by an orthonormal set of eigenvectors of the averaged periodogram matrix, where membership in the set is determined by a partitioning of the sorted observed eigenvalues into contiguous groups of sizes that match the dimensions of the corresponding true cointegrating subspaces. In Section 4, we show that the eigenvalues for the $k$th estimated cointegrating subspace are $O_p(n^{2d_k})$, where $n$ is the sample size and $d_k$ is the memory parameter of the cointegrating error for the $k$th true cointegrating subspace. This result and its refinements play a key role in our subsequent theory.

We will show in Theorem 1 that any vector in the $k$th estimated cointegrating subspace is, with high probability, close to the $k$th true cointegrating subspace, in the sense that the norm of the sine of the angle between these two subspaces converges in probability to zero. The norm of the sine of this angle is $O_p(n^{-\alpha_k})$, where $\alpha_k$ is the shortest distance between the memory parameters corresponding to the given and adjacent subspaces. This implies that the sine of the angle between any vector in the $k$th estimated cointegrating subspace and the $k$th true cointegrating subspace is $O_p(n^{-\alpha_k})$. (We provide more details on the notion of the sine of the angle between subspaces, and also the sine of the angle between a vector and a subspace, in Section 5.) This convergence rate, which improves as $\alpha_k$ increases, is at least as fast as the rates obtained for existing semiparametric estimators of cointegrating vectors in the bivariate case (see, e.g., [19] and the discussion in [3]), but not as fast as the parametric rate obtained by Hualde and Robinson [7] of $O_p(n^{-1/2})$ in the bivariate asymptotically stationary case if the difference $(\alpha_1)$ between the memory parameters of the observed series and the cointegrating error is less than $1/2$. Furthermore, we show in Lemma 15 that the normalized eigenvectors of the averaged periodogram matrix converge in distribution to random vectors that lie in the corresponding cointegrating subspace.

We then show in Section 6 that the cointegrating residuals corresponding to an estimated cointegrating vector can be used to obtain a consistent and asymptotically normal estimate of the memory parameter for the given cointegrating subspace, using a univariate Gaussian semiparametric estimator with a bandwidth that tends to $\infty$ more slowly than $n$. We also describe a procedure for consistently identifying the cointegrating subspaces, that is, for determining the number of subspaces and their dimensions. In Section 7, we provide a test for fractional cointegration which is appropriate for our model.

## 2. A fractional common components model.

Suppose that the original data are a $q \times 1$ time series such that the $(p-1)$th differences $\{y_t\}$ are weakly stationary with a common memory parameter $d_0 \in (-p+1/2, 1/2)$, where $p \geq 1$ is a fixed integer. The use of $(p-1)$th differences converts any additive polynomial trend of order $p-1$ in the original series into an additive constant. The value of this constant is irrelevant for our purposes since the estimators considered here are functions of the discrete Fourier transform at nonzero Fourier frequencies. We can, therefore, take the mean of $\{y_t\}$ to be zero, without loss of generality, and our estimators are invariant to polynomial trends of order $p-1$ in the original series.

In order to guarantee that the cointegrating relationships in the stochastic component of the levels are preserved in the differences, we apply a taper to the differences, that is, we multiply the differences by a sequence of constants prior to Fourier transformation. This prevents detrimental leakage effects due to potential overdifferencing and allows us to obtain uniform results over a wide range of memory parameters. A convenient family of tapers for use on the differences, and which we will use here, was given in Hurvich and Chen [8]. The exact form of the taper is given below.

The fractional common-components model for the $(q \times 1)$ series $\{y_t\}$ with cointegrating rank $r$ $(1 \le r < q)$ and $s$ cointegrating subspaces $(1 \le s \le r)$ is given by

$$y_t = \mathbf{A}_0 u_t^{(0)} + \mathbf{A}_1 u_t^{(1)} + \cdots + \mathbf{A}_s u_t^{(s)}, \tag{1}$$

where $\mathbf{A}_k$ $(0 \le k \le s)$ are $q \times a_k$ full-rank matrices with $a_0 = q - r$ and $a_1 + \cdots + a_s = r$ such that all columns of $\mathbf{A}_0, \ldots, \mathbf{A}_s$ are linearly independent, and $\{u_t^{(k)}\}$, $k = 0, \ldots, s$, are $(a_k \times 1)$ processes with memory parameters $\{d_k\}_{k=0}^{s}$ with $-p + 1/2 < d_s < \cdots < d_0 < 1/2$. Equation (1) can be written as

$$y_t = \mathbf{A} z_t, \tag{2}$$

where $z_t = \operatorname{vec}(u_t^{(0)}, \ldots, u_t^{(s)})$ and $\mathbf{A} = [\mathbf{A}_0 \ \ldots \ \mathbf{A}_s]$. We will make additional assumptions on $\{z_t\}$ in Section 3. These assumptions guarantee that $\{z_t\}$ is not cointegrated. The methodology presented in this paper does not require either $r$ or $s$ to be known.

Remark 1. Our assumption that all entries of $\{y_t\}$ have memory parameter $d$ implies that all rows of $A_0$ are nonzero. The model (1), without the assumption that all entries of $\{y_t\}$ have a common memory parameter, could also be entertained (though we do not pursue this here) and would then include the model considered by Robinson and Yajima [20].

Next, we exhibit the cointegrating subspaces. For any matrix $\mathbf{A}$, let $\mathcal{M}(\mathbf{A})$ denote the column space of $\mathbf{A}$ and let $\mathcal{M}^{\perp}(\mathbf{A})$ denote the orthogonal complement of $\mathbf{A}$. Note that for $k = 1, \ldots, s$,

$$\mathcal{M}^{\perp}(\mathbf{A}_0, \ldots, \mathbf{A}_k) \subset \mathcal{M}^{\perp}(\mathbf{A}_0, \ldots, \mathbf{A}_{k-1}).$$

Let $\mathcal{B}_0 = \mathcal{M}(\mathbf{A}_0)$ and $\mathcal{B}_k$, $k = 1, \ldots, s$, be the subspace such that

$$\mathcal{M}^{\perp}(\mathbf{A}_0, \ldots, \mathbf{A}_{k-1}) = \mathcal{M}^{\perp}(\mathbf{A}_0, \ldots, \mathbf{A}_k) \oplus \mathcal{B}_k$$

and $\mathcal{B}_k \perp \mathcal{M}^{\perp}(\mathbf{A}_0, \dots, \mathbf{A}_k)$. Hence, a nonzero vector $\beta \in \mathcal{B}_k$, $k \in \{1, \dots, s\}$, satisfies $\beta' \mathbf{A}_\ell = 0$, $\ell = 0, \dots, k-1$, and $\beta' \mathbf{A}_k \neq 0$. Also, $\mathcal{B}_j \perp \mathcal{B}_k$ for $j \neq k$, $(j,k) \in \{0, \dots, s\}$ and

$$\mathbb{R}^q = \mathcal{B}_0 \oplus \mathcal{B}_1 \oplus \cdots \oplus \mathcal{B}_s. \tag{3}$$

It can be seen from (1) and the preceding discussion that any nonzero vector $\beta \in \mathcal{B}_k$ with $k \in \{1, \dots, s\}$ produces a cointegrating error series $\{\beta' y_t\}$ with memory parameter $d_k$. Thus, $\mathcal{B}_1, \dots, \mathcal{B}_s$ are the cointegrating subspaces. The space $\mathcal{B}_0$, on the other hand, is the space spanned by any basis of non-cointegrating vectors in $\mathbb{R}^q$. Equation (3) shows that $\mathbb{R}^q$ may be written as a direct sum of the space of non-cointegrating vectors and the space of cointegrating vectors, and that the latter space may be further decomposed into a direct sum of cointegrating subspaces.

**3. Assumptions.** Here, we specify a linear model for the series $z_t = \text{vec}(u_t^{(0)}, \dots, u_t^{(s)})$. As stated in the previous section, we assume that $\{u_t^{(k)}\}$, $k = 0, \dots, s$, are $(a_k \times 1)$ processes with memory parameters $\{d_k\}_{k=0}^{s}$ with $-p + 1/2 < d_s < \cdots < d_0 < 1/2$. Define $N_0 = \{1, \dots, a_0\}$ and $N_k = \{(a_0 + \cdots + a_{k-1}) + 1, \dots, (a_0 + \cdots + a_k)\}$ for $k = 1, \dots, s$. Our results in this paper assume $s > 0$, unless explicitly stated otherwise.

Let $\psi_k$ be a sequence of $q \times q$ matrices such that

$$\psi_k = \frac{1}{2\pi} \int_{-\pi}^{\pi} e^{ik\omega} \mathbf{\Psi}(\omega)\, d\omega,$$

where for each $\omega \in [-\pi, \pi]$, $\mathbf{\Psi}(\omega)$ is a complex-valued matrix such that $\mathbf{\Psi}(-\omega) = \overline{\mathbf{\Psi}(\omega)}$ and where $\psi_0$ is an identity matrix.

Define the $q \times 1$ vector process $\{z_t\}$ as

$$z_t = \sum_{k=-\infty}^{\infty} \psi_k \varepsilon_{t-k}, \tag{4}$$

where $\{\varepsilon_t = (\varepsilon_{t,1}, \dots, \varepsilon_{t,q})'\} \sim iid(\mathbf{0}, 2\pi\mathbf{\Sigma})$, $\mathbf{\Sigma}$ is a symmetric positive definite matrix with entries $\sigma_{ab}$, $a, b \in \{1, \dots, q\}$, and $E\|\varepsilon_t\|^4 < \infty$, where $\|\cdot\|$ denotes the Euclidean norm. The spectral density matrix of $\{z_t\}$ is

$$\mathbf{f}(\omega) = \mathbf{\Psi}(\omega)\mathbf{\Sigma}\mathbf{\Psi}^*(\omega), \qquad \omega \in [-\pi, \pi],$$

where the superscript $*$ denotes conjugate transposition. We further assume that for $\omega \in [-\pi, \pi]$, the $(a,b)$th entry of $\mathbf{\Psi}(\omega)$ is given by

$$\Psi_{ab}(\omega) = (1 - e^{-i\omega})^{-d_{ab}} \tau_{ab}(\omega) e^{i\phi_{ab}(\omega)}, \tag{5}$$

where $d_{aa} = d_k$ for $a \in N_k$, $d_{ab} \leq \min(d_k, d_h)$ for $a \in N_k$, $b \in N_h$, $b \neq a$, $k, h = 0, \dots, s$, and for all $a, b \in \{1, \dots, q\}$, $\tau_{ab}(\cdot)$ are positive even real-valued functions and $\phi_{ab}(\cdot)$ are odd real-valued functions, all continuously differentiable in an interval containing zero. It follows from (5) that the first derivatives of $\Psi_{ab}(\omega)$ satisfy

$$\Psi'_{ab}(\omega) = O(|\Psi_{aa}(\omega)\Psi_{bb}(\omega)|^{1/2}|\omega|^{-1}). \tag{6}$$

In keeping with (5), we assume that we can write the spectral density matrix of $\{z_t\}$ as

$$\mathbf{f}(\omega) = \boldsymbol{\Upsilon}(\omega)\mathbf{f}^{\dagger}(\omega)\boldsymbol{\Upsilon}^{*}(\omega), \tag{7}$$

where $\boldsymbol{\Upsilon}(\omega) = \operatorname{diag}\{(1-e^{-i\omega})^{-d_0}, \ldots, (1-e^{-i\omega})^{-d_0}, \ldots, (1-e^{-i\omega})^{-d_s}, \ldots, (1-e^{-i\omega})^{-d_s}\}$, that is, the $a$th diagonal entry is $(1-e^{-i\omega})^{-d_k}$ for all $a \in N_k$, $k = 0, \ldots, s$, and

$$\mathbf{f}^{\dagger}(\omega) = \boldsymbol{\Psi}^{\dagger *}(\omega)\Sigma\boldsymbol{\Psi}^{\dagger}(\omega) \tag{8}$$

is positive definite, Hermitian, continuous at zero frequency and, therefore, real-valued at zero frequency with $\Psi^{\dagger}_{ab}(\omega) = \tau_{ab}(\omega)e^{i\phi_{ab}(\omega)}$. Thus, $\{z_t\}$ is not fractionally cointegrated (see [18]).

**4. The averaged periodogram matrix and its eigenvalues.** For any vector sequence of observations $\{\xi_t\}_{t=1}^{n}$, define the tapered discrete Fourier transform by

$$J_{\xi}(\omega_j) = \frac{1}{\sqrt{2\pi\sum|h_t^{p-1}|^2}}\sum_{t=1}^{n} h_t^{p-1}\xi_t e^{i\omega_j t},$$

where $\omega_j = 2\pi j/n$ is the $j$th Fourier frequency and $\{h_t\}$ is the complex-valued taper of Hurvich and Chen [8],

$$h_t = 0.5(1 - e^{i2\pi t/n}), \qquad t = 1, \ldots, n.$$

Note that $p = 1$ yields the no-tapering case. Next, define the tapered cross-periodogram matrix of two vector sequences $\{\xi_t\}_{t=1}^{n}$ and $\{\zeta_t\}_{t=1}^{n}$ by

$$I_{\xi\zeta}(\omega_j) = J_{\xi}(\omega_j)J_{\zeta}^{*}(\omega_j).$$

We will work with the (real part of the) averaged periodogram matrix of a sample of $n$ observations $\{y_t\}_{t=1}^{n}$,

$$I_m = \sum_{j=1}^{m}\operatorname{Re}\{I_{yy}(\omega_j)\},$$

where $m$ is a fixed positive integer, $m > q + 3$. (This condition is motivated in the proof of Lemma 8.)

Define $I_m(\xi_t,\zeta_t)=\sum_{j=1}^m \operatorname{Re}\{I_{\xi\zeta}(\omega_j)\}$. We first focus on the asymptotic distribution of $I_m(z_t,z_t)$. Define the function (for $x\in\mathbb{R}$)

$$\Delta_p(x)=\binom{2p-2}{p-1}^{-1/2}\sum_{k=0}^{p-1}\binom{p-1}{k}(-1)^k\Delta(x+2\pi k),$$

where

$$\Delta(x)=\frac{1}{\sqrt{2\pi}}\frac{e^{ix}-1}{ix}.$$

Now, define

$$\upsilon_j(x)=\frac{1}{2}[\overline{\Delta_p(-x+2\pi j)}+\Delta_p(x+2\pi j)],$$
$$\nu_j(x)=\frac{i}{2}[\overline{\Delta_p(-x+2\pi j)}-\Delta_p(x+2\pi j)].$$

Define the Hermitian positive definite $q\times q$ matrix-valued measure $\mathbf{G}_0$ on $\mathbb{R}$ by

$$\mathbf{G}_0(dx)=\mathbf{\Pi}(x)\mathbf{f}^{\dagger}(0)\mathbf{\Pi}^*(x)\,dx \tag{9}$$

for $x>0$ and $\mathbf{G}_0(-dx)=\overline{\mathbf{G}_0(dx)}$, where

$$\mathbf{\Pi}(x)=\operatorname{diag}(e^{-i\pi d_0/2}|x|^{-d_0},\dots,e^{-i\pi d_0/2}|x|^{-d_0},\dots,e^{-i\pi d_s/2}|x|^{-d_s},\dots,\\ e^{-i\pi d_s/2}|x|^{-d_s}).$$

Let $\mathbf{U}_n$ and $\mathbf{V}_n$ be $q\times m$ matrices given by

$$\mathbf{U}_n=\mathbf{d}_n^{-1}\operatorname{Re}(J_{z,1},\dots,J_{z,m})\quad\text{and}\quad\mathbf{V}_n=\mathbf{d}_n^{-1}\operatorname{Im}(J_{z,1},\dots,J_{z,m}). \tag{10}$$

LEMMA 1. *Let $\mathbf{d}_n$ be a $(q\times q)$ diagonal matrix with ith diagonal entry $n^{d_k}$, $i\in N_k$, $k=0,\dots,s$ and $\mathbf{Q}_n=\mathbf{d}_n^{-1}I_m(z_t,z_t)\mathbf{d}_n^{-1}=(\mathbf{U}_n,\mathbf{V}_n)(\mathbf{U}_n,\mathbf{V}_n)'$. If $m\ge q$, then*

$$\mathbf{Q}_n\xrightarrow{D}\mathbf{UU}'+\mathbf{VV}',$$

*where $\mathbf{U}=(U_1,\dots,U_m)$ and $\mathbf{V}=(V_1,\dots,V_m)$, $U_j$, $V_k$ are $q\times 1$ vectors and* $\operatorname{vec}(\mathbf{U},\mathbf{V})$ *is a $2mq$-variate normal random variable with zero mean and covariance matrix $\mathbf{\Xi}$ determined by*

$$\mathbb{E}(U_jU_k')=\int_{\mathbb{R}}\upsilon_j(x)\overline{\upsilon_k(x)}\mathbf{G}_0(dx),$$
$$\mathbb{E}(V_jV_k')=\int_{\mathbb{R}}\nu_j(x)\overline{\nu_k(x)}\mathbf{G}_0(dx),$$
$$\mathbb{E}(U_jV_k')=\int_{\mathbb{R}}\upsilon_j(x)\overline{\nu_k(x)}\mathbf{G}_0(dx).$$

*Furthermore, $\mathbf{UU}'+\mathbf{VV}'$ is positive definite and has distinct eigenvalues with probability 1.*

PROOF. The proof is identical to the proof of Lemma 1, Corollary 1 and 2 of [4]. □

We next derive upper and lower bounds for the eigenvalues of $I_m(y_t, y_t)$. We will use the notation $\lambda_j(\cdot)$ for the $j$th eigenvalue of a given Hermitian matrix, $\lambda_j(\cdot) \geq \lambda_{j+1}(\cdot)$. Also, we let $\lambda_j = \lambda_j(I_m(y_t, y_t))$. We have the following lemma:

LEMMA 2. *$\lambda_j = O_p(n^{2d_k})$ for $j \in N_k$, $k = 0, \dots, s$.*

In the case $k \geq 1$, the upper bound in Lemma 2 strengthens Lemma 4 of [4].

LEMMA 3. *Let $j_k^* = \max\{j : j \in N_k\}$ and let $\mathbf{Q}_n^{(k)}$ be the leading $j_k^* \times j_k^*$ principal submatrix of $\mathbf{Q}_n$ for $k = 0, \dots, s$. Then*

$$n^{-2d_k} \lambda_{j_k^*} \geq c_k \lambda_{j_k^*}(\mathbf{Q}_n^{(k)}) \xrightarrow{D} \eta_{j_k^*}^{(k)},$$

*where $c_k > 0$ and $\eta_{j_k^*}^{(k)}$ is a random variable that has no mass at 0.*

**5. Estimation of the cointegrating subspaces.** Let $\mathbf{X}(\cdot) = [\chi_1(\cdot) \ \dots \ \chi_q(\cdot)]$ be an orthogonal matrix such that $\chi_j(\cdot)$ is the eigenvector corresponding to the $j$th largest eigenvalue $\lambda_j(\cdot)$ of a given symmetric $q \times q$ matrix and let $\mathbf{X}_k(\cdot)$ be a matrix with columns $\chi_j(\cdot)$, $j \in N_k$, for $k = 0, \dots, s$. Also, let $\chi_j = \chi_j(I_m(y_t, y_t))$, $\mathbf{X} = \mathbf{X}(I_m(y_t, y_t))$ and $\mathbf{X}_k = \mathbf{X}_k(I_m(y_t, y_t))$. For $k = 0, 1, \dots, s$, let $\mathbf{B}_k$ be a $q \times a_k$ matrix with orthonormal columns such that $\mathcal{M}(\mathbf{B}_k) = \mathcal{B}_k$ and let $\mathbf{B} = [\mathbf{B}_0 \ \dots \ \mathbf{B}_s]$. Since $\mathbf{B}'\mathbf{B} = \mathbf{I}$, it follows that for any $q \times q$ matrix $\mathbf{P}$, $\mathbf{B}'\mathbf{P}\mathbf{B}$ is similar to $\mathbf{P}$, that is, $\lambda_j(\mathbf{P}) = \lambda_j(\mathbf{B}'\mathbf{P}\mathbf{B})$ and $\chi_j(\mathbf{P}) = \mathbf{B}'\chi_j(\mathbf{B}'\mathbf{P}\mathbf{B})$.

Define

$$\mathbf{\Phi} = \mathbf{B}' I_m(y_t, y_t) \mathbf{B}$$

and partition $\mathbf{\Phi}$ into $(s+1)^2$ blocks such that the $(k, \ell)$ block $\mathbf{\Phi}_{k\ell}$ has dimension $(a_k \times a_\ell)$ for $k, \ell = 0, \dots, s$. Define $\mathbf{\Phi}_D = \operatorname{diag}[\mathbf{\Phi}_{00}, \dots, \mathbf{\Phi}_{ss}]$ and $\Delta\mathbf{\Phi} = \mathbf{\Phi} - \mathbf{\Phi}_D$, so that

$$\mathbf{\Phi} = \mathbf{\Phi}_D + \Delta\mathbf{\Phi}.$$

We have

$$I_m(y_t, y_t) = \mathbf{B}\mathbf{\Phi}\mathbf{B}' = \mathbf{B}\mathbf{\Phi}_D\mathbf{B}' + \mathbf{B}\Delta\mathbf{\Phi}\mathbf{B}' =: \mathbf{H} + \Delta\mathbf{H},$$

so we can think of $I_m(y_t, y_t)$ as a perturbed version of $\mathbf{H}$. Using results of Barlow and Slapničar [2] on perturbation theory for eigenvalues and eigenvectors of nonrandom Hermitian matrices, we will show in Lemma 4 that the $k$th estimated cointegrating subspace $\mathcal{M}(\mathbf{X}_k)$ is close to $\mathcal{M}(\mathbf{X}_k(\mathbf{H}))$ in the sense that the norm of the sine of the angle between the two subspaces converges to 0 in probability.

Let $\Theta(\cdot,\cdot)$ denote the matrix of canonical angles between two subspaces of the same dimension (see, e.g., [22], page 43). The notion of the sine of the angle between two subspaces of the same dimension is given in [5]. For simplicity, suppose that $\mathbf{S}$ and $\mathbf{T}$ are both real $q \times a$ matrices $(q > a)$ with orthonormal columns. Then the orthogonal projector into $\mathcal{M}(\mathbf{T})$ is given by $\mathbf{T}\mathbf{T}'$ and the projector into the orthogonal complement $\mathcal{M}^{\perp}(\mathbf{T})$ of $\mathcal{M}(\mathbf{T})$ is given by $\mathbf{I} - \mathbf{T}\mathbf{T}'$, where $\mathbf{I}$ is a $q \times q$ identity matrix. The sine of the angle between $\mathcal{M}(\mathbf{S})$ and $\mathcal{M}(\mathbf{T})$ is an $a \times a$ matrix defined in [5] and denoted by $\sin\Theta(\mathcal{M}(\mathbf{S}),\mathcal{M}(\mathbf{T}))$. It follows from [5], page 10 that $\|\sin\Theta(\mathcal{M}(\mathbf{S}),\mathcal{M}(\mathbf{T}))\|_F = \|(\mathbf{I}-\mathbf{T}\mathbf{T}')\mathbf{S}\mathbf{S}'\|_F$, where $\|\cdot\|_F$ is the Frobenius norm. It follows from [22], Corollary 5.4, page 43 that

$$\|\sin\Theta(\mathcal{M}(\mathbf{S}),\mathcal{M}(\mathbf{T}))\|_F = \|(\mathbf{T}^{\perp})'\mathbf{S}\|_F, \tag{11}$$

where $\mathbf{T}^{\perp}$ is a matrix with orthonormal columns spanning $\mathcal{M}^{\perp}(\mathbf{T})$, so that $\|(\mathbf{T}^{\perp})'\mathbf{S}\|_F$ is the square root of the sum of the squared lengths of the residuals from the orthogonal projections of the columns of $\mathbf{S}$ on the space $\mathcal{M}(\mathbf{T})$.

For any nonzero vector $x \in \mathcal{M}(\mathbf{S})$, the sine of the angle between $x$ and the subspace $\mathcal{M}(\mathbf{T})$ is a real number defined as

$$\sin\theta(x,\mathcal{M}(\mathbf{T})) = \frac{\|(\mathbf{I}-\mathbf{T}\mathbf{T}')x\|}{\|x\|},$$

see [24], page 274. It then follows from (11) that

$$\max_{x\in\mathcal{M}(\mathbf{S})} |\sin\theta(x,\mathcal{M}(\mathbf{T}))| \le \|(\mathbf{T}^{\perp})'\mathbf{S}\|_F.$$

In Lemma 5, we show that under the additional assumption that the process is Gaussian, $\mathcal{M}(\mathbf{X}_k(\mathbf{H}))$ is equal to $\mathcal{B}_k$ with probability approaching one, for $k = 0,\dots,s$. Lemmas 4 and 5, taken together, imply our Theorem 1, stating that if the process is Gaussian, then the $k$th estimated cointegrating subspace $\mathcal{M}(\mathbf{X}_k)$ is close to the corresponding true cointegrating subspace $\mathcal{B}_k$, in the sense that $\|\sin\Theta\{\mathcal{M}(\mathbf{X}_k),\mathcal{B}_k\}\|_F = O_p(n^{-\alpha_k})$, where $\alpha_k$ is the shortest distance between the memory parameters corresponding to the given and adjacent subspaces, that is,

$$\alpha_k = \begin{cases} d_0 - d_1, & k = 0,\\ \min\{(d_{k-1}-d_k),(d_k-d_{k+1})\}, & k = 1,\dots,s-1,\\ d_{s-1}-d_s, & k = s.\end{cases}$$

Lemma 4. *The sine of the angle between $\mathcal{M}(\mathbf{X}_k)$ and $\mathcal{M}(\mathbf{X}_k(\mathbf{H}))$ satisfies*

$$\|\sin\Theta\{\mathcal{M}(\mathbf{X}_k(\mathbf{H})),\mathcal{M}(\mathbf{X}_k)\}\|_F = O_p(n^{-\alpha_k}).$$

The following Gaussianity assumption is sufficient for obtaining a rate at which $P(\mathcal{M}(\mathbf{X}_k(\mathbf{H})) \neq \mathcal{B}_k)$ converges to zero. More specifically, the assumption allows us to bound the inverse second moment of eigenvalues of $\mathbf{Q}_n$. We believe that such bounds, and therefore Lemma 5, hold without the Gaussianity assumption, but we will not pursue this here.

ASSUMPTION 1. *The process* $\{\varepsilon_t\}$ *in* (4) *is Gaussian.*

LEMMA 5. *Under Assumption* 1, $P(\mathcal{M}(\mathbf{X}_k(\mathbf{H})) \neq \mathcal{B}_k) = O(n^{-2\alpha_k})$, $k = 0, \dots, s$.

The following theorem is a corollary of Lemmas 4 and 5:

THEOREM 1. *Under Assumption* 1,

$$\|\sin\Theta\{\mathcal{M}(\mathbf{X}_k), \mathcal{B}_k\}\|_F = O_p(n^{-\alpha_k}), \qquad k = 0, \dots, s.$$

**6. Estimation of the memory parameters using cointegrating residuals.** Let $b = \chi_a$, where $a \in \{1, \dots, q\}$. Recall that $\chi_a$ is the eigenvector of $I_m(y_t, y_t)$ corresponding to the $a$th largest eigenvalue of the matrix. Then there exists a uniquely defined value $k \in \{0, \dots, s\}$ such that $a \in N_k$. Note that $k$ is fixed but unknown. We then use this vector $b$ to construct the residual process $\{v_t\}$, where

$$v_t := b' y_t = b'\mathbf{A}_0 u_t^{(0)} + b'\mathbf{A}_1 u_t^{(1)} + \cdots + b'\mathbf{A}_k u_t^{(k)} + \cdots + b'\mathbf{A}_s u_t^{(s)}. \tag{12}$$

The periodogram of $\{v_t\}$ is $I_{vv}(\omega_j) = b'\mathbf{A} I_{zz}(\omega_j)\mathbf{A}' b$. We consider the Gaussian semiparametric estimator or GSE (see [10, 17]) for $d_{aa}$ [see (5)] based on $\{v_t\}$,

$$\hat{d}_{aa} = \arg\min_{d \in \Theta} R(d) = \arg\min_{d \in \Theta} \left( \log \widehat{G}(d) - 2d \left( \frac{1}{m_n} \sum_{j=1}^{m_n} \log \omega_{\tilde{j}} \right) \right), \tag{13}$$

where $\Theta = [\Delta_1, \Delta_2]$, $-p + 0.5 < \Delta_1 < \Delta_2 < 0.5$, $\omega_{\tilde{j}} = 2\pi\tilde{j}/n$, $\tilde{j} = j + (p-1)/2$ and

$$\widehat{G}(d) = \frac{1}{m_n} \sum_{j=1}^{m_n} \frac{I_{vv}(\omega_j)}{\omega_{\tilde{j}}^{-2d}} = \frac{1}{m_n} \sum_{j=1}^{m_n} \frac{b'\mathbf{A} I_{zz}(\omega_j)\mathbf{A}' b}{\omega_{\tilde{j}}^{-2d}}.$$

Here, we use slightly shifted Fourier frequencies $\omega_{\tilde{j}}$ to parallel corresponding shifts inherent in our tapering scheme and thereby reduce finite-sample bias, as was also done in [8].

The two theorems below establish the consistency and the limiting distribution of the $\hat{d}_{aa}$, under some additional conditions on the transfer function $\Psi^{\dagger}_{ab}(\omega) = \tau_{ab}(\omega) e^{i\phi_{ab}(\omega)}$; see (5). Following [9], we define a smoothness class for transfer functions as follows. For $\mu > 1$ and $1 < \rho \leq 2$, let $\mathcal{L}^*(\mu, \rho)$ be the set of continuously differentiable functions $u$ on $[-\pi, \pi]$ such that for all $x, y$ with $|x| \in (0, \pi]$, $|y| \in (0, \pi]$,

$$\frac{\max_{0\leq z\leq\pi}|u(z)|}{\min_{0\leq z\leq\pi}|u(z)|} \leq \mu, \qquad \frac{|u(x)-u(y)|}{\min_{0\leq z\leq\pi}|u(z)|} \leq \mu\frac{|y-x|}{\min(|x|,|y|)},$$

$$\frac{|u'(x)-u'(y)|}{\min_{0\leq z\leq\pi}|u(z)|} \leq \mu\frac{|y-x|^{(\rho-1)}}{[\min(|x|,|y|)]^{\rho}}.$$

It follows from the discussion in [9] that if $\Psi^{\dagger}_{aa}(\omega)$ is the transfer function of a stationary and invertible autoregressive moving average process, or of the short-memory component of a stationary and invertible fractional Gaussian noise with a suitable choice of the moving average representation, then $\Psi^{\dagger}_{aa}(\omega) \in \mathcal{L}^*(\mu, \rho)$ for some $\mu$, with $\rho = 2$.

We now state an assumption on $\mathbf{\Psi}^{\dagger}$.

Assumption 2. *For all $a, b \in \{1, \ldots, q\}$, $\Psi^{\dagger}_{ab} \in \mathcal{L}^*(\mu, \rho)$ for some $\mu > 1$ and some $\rho \in (1, 2]$.*

Note that this assumption is global in that it pertains to the behavior of $\mathbf{\Psi}^{\dagger}$ at all frequencies. By contrast, our estimation of the $d_{aa}$ is based on frequencies in a shrinking neighborhood around zero. It seems plausible, then, that a local version of Assumption 2 would suffice for our purposes, although we do not pursue this here.

The following standard assumption is needed to establish the consistency of $\hat{d}_{aa}$:

Assumption 3A. *As $n \to \infty$,*

$$\frac{1}{m_n} + \frac{m_n}{n} \to 0.$$

Theorem 2. *Under Assumptions 1, 2 and 3A, for $a \in \{1, \ldots, q\}$, $\hat{d}_{aa} \xrightarrow{p} d_{aa}$.*

The next assumption is used for establishing the asymptotic normality of $m_n^{1/2}(\hat{d}_{aa} - d_{aa})$.

Assumption 3B. *Suppose that $a \in N_k$.*

(i) *If $k \in \{1, \ldots, s\}$, then $d_{k-1} - d_k > 1/2$.*

(ii) *If* $k \in \{0, \ldots, s-1\}$, *then as* $n \to \infty$,

$$\frac{1}{m_n} + \frac{m_n^{1+2(d_k - d_{k+1})} \log^2 m_n}{n^{2(d_k - d_{k+1})}} \to 0.$$

Note that part (i) is vacuous if $k = 0$ and part (ii) is vacuous if $k = s$. Assumption 3B may be compared with the assumptions in Theorems 2 and 4 of Velasco [23] which the author required for residual-based estimators of the memory parameters of a bivariate fractionally cointegrated system. The problem here is that a linear combination of series with slightly different memory parameters will typically have an irregular short-memory component in its spectral density.

To present the asymptotic variance of $\hat{d}_{aa}$, we define

$$\Phi_p = \frac{\Gamma(4p-3)\Gamma^4(p)}{\Gamma^4(2p-1)}.$$

THEOREM 3. *Under Assumptions* 1, 2 *and* 3B, *for* $a \in \{1, \ldots, q\}$,

$$m_n^{1/2}(\hat{d}_{aa} - d_{aa}) \xrightarrow{D} N(0, \Phi_p/4).$$

Note that in Theorem 3, the limiting distribution of $m_n^{1/2}(\hat{d}_{aa} - d_{aa})$ has mean zero. This asymptotic unbiasedness is ensured by Assumption 3B, which places strong restrictions on the separation between the memory parameters and also places a potentially stringent upper bound on the bandwidth $m_n$. A much weaker and, indeed, more standard assumption involving only $m_n$ is the following:

ASSUMPTION 3C. *As* $n \to \infty$,

$$\frac{1}{m_n} + \frac{m_n^{1+2\rho} \log^2 m_n}{n^{2\rho}} \to 0.$$

If we account for the asymptotic bias, which can be determined from Lemma 20, and use Assumption 3C, we obtain the following result:

COROLLARY 1. *Suppose* $a \in N_k$, *where* $k \in \{0, \ldots, s\}$. *Under Assumptions* 1, 2 *and* 3C, *we have*

$$m_n^{1/2}(\hat{d}_{aa} - d_{aa} - \mu_n) \xrightarrow{D} N(0, \Phi_p/4),$$

*where* $\mu_n = O_p(m_n^{d_k - d_{k-1}} + \omega_{m_n}^{d_k - d_{k+1}})$, *the* $O_p(m_n^{d_k - d_{k-1}})$ *term is vacuous if* $k = 0$ *and the* $O_p(\omega_{m_n}^{d_k - d_{k+1}})$ *term is vacuous if* $k = s$.

Here, we present some results on the vector of GSE-estimated memory parameters, $\hat{d} = (\hat{d}_{11}, \ldots, \hat{d}_{qq})'$, which is an estimate of $d = (d_{11}, \ldots, d_{qq})'$. Let $w_t = \mathbf{X}' y_t$ be the $q \times 1$ residual vector so that the entries of $\hat{d}$ are based on those of $w_t$. Note that by Lemma 15, $\mathbf{X} \xrightarrow{D} \mathring{\mathbf{X}}$, where $\mathring{\mathbf{X}}$ is a continuous function of $\mathbf{U}$ and $\mathbf{V}$ in Lemma 1. We will need the following assumption for our results:

ASSUMPTION 3D. (i) *For all* $k \in \{0, \ldots, s\}$, $\alpha_k > 1/2$.
(ii) *As* $n \to \infty$,

$$\frac{1}{m_n} + \frac{m_n^{1+2\xi} \log^2 m_n}{n^{2\xi}} \to 0,$$

*where* $\xi = \min\{\min_k \alpha_k, \rho\}$.

COROLLARY 2. *Under Assumptions 1, 2 and 3D,*

$$m_n^{1/2}(\hat{d} - d) \xrightarrow{D} N\left(0, \frac{\Phi_p}{4} (diag\, \mathbf{\Omega})^{-1} \circ \mathbf{\Omega} \circ \mathbf{\Omega} \circ (diag\, \mathbf{\Omega})^{-1}\right),$$

*where*

$$\mathbf{\Omega} = \mathbb{E}(\mathring{\mathbf{X}}' \mathbf{A} \mathbf{f}^{\dagger}(0) \mathbf{A}' \mathring{\mathbf{X}}).$$

REMARK 2. Simulation results not shown here reveal that the small-sample bias is reduced and the variance is stabilized if the GSE estimators omit the first $m + p - 1$ frequencies. This does not affect the validity of Corollary 2. Note that if no frequencies are omitted, then the first $m+p-1$ frequencies are used twice: once for estimating the cointegrating vector and once for estimating the memory parameter. If the frequencies are omitted, the finite-sample approximation to the variance in [8] is quite accurate.

Corollary 2 yields the following result on the asymptotic distribution of $m_n^{1/2}(\hat{d}_{aa} - \hat{d}_{bb} - (d_{aa} - d_{bb}))$, under conditions that ensure asymptotic unbiasedness:

COROLLARY 3. *Under the assumptions of Corollary 2, for* $a, b \in \{1, \ldots, q\}$, $a \neq b$,

$$m_n^{1/2}(\hat{d}_{aa} - \hat{d}_{bb} - (d_{aa} - d_{bb})) \xrightarrow{D} N\left(0, \frac{\Phi_p}{2}\left(1 - \frac{\Omega_{ab}^2}{\Omega_{aa}\Omega_{bb}}\right)\right).$$

Next, we modify Corollary 3 to include a bias term, thereby allowing for weaker assumptions.

COROLLARY 4. *If $a \in N_k$, $b \in N_h$, for $k, h \in \{0, \dots, s\}$, then under the assumptions of Corollary 1,*

$$m_n^{1/2}(\hat{d}_{aa} - \hat{d}_{bb} - (d_{aa} - d_{bb}) - \tilde{\mu}_n) \xrightarrow{D} N\left(0, \frac{\Phi_p}{2}\left(1 - \frac{\Omega_{ab}^2}{\Omega_{aa}\Omega_{bb}}\right)\right),$$

*where*

$$\tilde{\mu}_n = O_p(m_n^{d_k - d_{k-1}} + m_n^{d_h - d_{h-1}} + \omega_{m_n}^{d_k - d_{k+1}} + \omega_{m_n}^{d_h - d_{h+1}}).$$

Given data from model (1), assumed to possess fractional cointegration, the number $s > 0$ of cointegrating subspaces and their dimensions $a_1, \dots, a_s$, as well as the dimension $a_0$ of the non-cointegrating space, will be unknown in general. Here, we assume Gaussianity. Let $\delta^* > 0$ be the minimum separation between the memory parameters, $\delta^* = \min(d_0 - d_1, \dots, d_{s-1} - d_s)$, and assume first that $\delta^* > 1/2$. We can compare the GSE estimators $\hat{d}_{jj}$ and $\hat{d}_{j+1,j+1}$ for $j = 1, \dots, q$ using a bandwidth $m_n$ satisfying Assumption 3D, part (ii), with $\xi = \min\{\delta^*, \rho\}$. Fix an $\varepsilon \in (0, 1/2)$ and a $C > 0$. Then, for each $j \in \{1, \dots, q-1\}$, we declare that $d_{jj} - d_{j+1,j+1} \neq 0$ if and only if $\hat{d}_{jj} - \hat{d}_{j+1,j+1} > C m_n^{-1/2+\varepsilon}$. This leads to a procedure for consistently identifying $s$, $a_0, \dots, a_s$, which can be justified by Corollary 3. A more complicated identification procedure, justified by Corollary 4, may be constructed if $\delta^* \leq 1/2$. This requires further tuning parameters which depend on $\delta^*$, owing to the presence of the nonstandard term $\tilde{\mu}_n$, which increases as the separation of the relevant memory parameters decreases. Unfortunately, such procedures are of limited practical value as $\delta^*$ will be unknown in practice. Nevertheless, we note that lower bounds on $\delta^*$ (assuming $s = 1$) arise implicitly or explicitly in other works on semiparametric fractional cointegration. (See [20], Assumption D and [23], Theorems 2 and 4.)

**7. Testing for fractional cointegration.** In model (1), used throughout the paper thus far, we have assumed that $s \geq 1$ so that cointegration exists. Here, we expand model (1) to include the case of no cointegration ($s = 0$, or equivalently, $r = 0$), that is,

$$y_t = \mathbf{A}_0 u_t^{(0)}, \tag{14}$$

where $\mathbf{A}_0$ is $q \times q$ with linearly independent columns and all entries of $u_t^{(0)}$ have memory parameter $d_0$.

In practice, it is of interest to test for the presence of fractional cointegration. Such a test was proposed by Marinucci and Robinson ([13], pages 236–237), following on from an idea originally suggested in a different context by Hausman [6], using a comparison of two estimates of $d_0$, one based on a multivariate Gaussian semiparametric estimator (see [11]) using $\{y_t\}_{t=1}^n$ with an imposed restriction that all entries have the same memory parameter, and the other estimator based on a univariate Gaussian semiparametric estimator of $d_0$ using (say) the first entry $\{y_{1,t}\}$ of $\{y_t\}$. It seems possible to use this idea, together with differencing and tapering, to yield a test for fractional integration in the current context, although we do not pursue this here. We focus instead on residual-based methods in which estimated memory parameters based on the various cointegrating residual series are compared. In a bivariate context, Velasco [23] has considered properties of semiparametric memory parameter estimates based on cointegrating residuals under certain assumptions on the rate of convergence of the semiparametric estimator of the cointegrating parameters. However, the author did not present a test for cointegration since his assumptions ruled out the no-cointegration case.

For our GSE estimators $\hat{d}$ based on cointegrating residuals, we have the following extensions of Corollary 2 and Corollary 3 to the no-cointegration case (14):

THEOREM 4. *Under Assumptions 1, 2 and 3C, if there is no cointegration, then*

$$m_n^{1/2}(\hat{d}-d) \xrightarrow{D} N\left(0, \frac{\Phi_p}{4}(diag\,\mathbf{\Omega})^{-1} \circ \mathbf{\Omega} \circ \mathbf{\Omega} \circ (diag\,\mathbf{\Omega})^{-1}\right).$$

COROLLARY 5. *Under Assumptions 1, 2 and 3C, if there is no cointegration, then for $a, b \in \{1, \ldots, q\}$,*

$$m_n^{1/2}(\hat{d}_{aa} - \hat{d}_{bb}) \xrightarrow{D} N\left(0, \frac{\Phi_p}{2}\left(1 - \frac{\Omega_{ab}^2}{\Omega_{aa}\Omega_{bb}}\right)\right).$$

Corollaries 4 and 5 justify a conservative hypothesis test for the null hypothesis of no cointegration based on the test statistic $T_n = m_n^{1/2}(\hat{d}_{11} - \hat{d}_{qq})$ whereby, for a nominal level $\alpha$ test, the null hypothesis is rejected in favor of the cointegration alternative hypothesis if and only if $T_n > (\Phi_p/2)^{1/2} z_{\alpha/2}$. Here, a bandwidth $m_n$ satisfying Assumption 3C should be used. The test is conservative since $(\Phi_p/2)$ is an upper bound for the asymptotic variance of $T_n$.

## 8. Proofs.

### 8.1. *Proofs for Section 4.*

PROOF OF LEMMA 2. Note that

$$\mathbf{\Phi} = \mathbf{B}'\mathbf{A} I_m(z_t, z_t)\mathbf{A}'\mathbf{B},$$

where $\mathbf{B}'\mathbf{A}$ is an upper triangular block matrix. We first partition $\mathbf{\Phi}$ into $(s+1)\times(s+1)$ blocks such that the $(k,\ell)$ block has dimension $(a_k \times a_\ell)$. Let $z_t^{(k)} = (u_t^{(k)}, \ldots, u_t^{(s)})$ and $\mathbf{A}^{(k)} = [\mathbf{A}_k \;\ldots\; \mathbf{A}_s]$, $k = 0, 1, \ldots, s$. We have

$$(15)\quad \begin{aligned} &\mathbf{\Phi}_{k\ell} = \mathbf{B}_k'\mathbf{A}^{(k)} I_m(z_t^{(k)}, z_t^{(\ell)})\mathbf{A}^{(\ell)\prime}\mathbf{B}_\ell, \qquad \text{for } k \le \ell,\; k,\ell = 0,1,\ldots,s,\\ &\mathbf{\Phi}_{\ell k} = \mathbf{\Phi}_{k\ell}'. \end{aligned}$$

Fix a value of $k \in \{0, \ldots, s\}$. Note that by Lemma 1, all the elements in the $k$th block, $\mathbf{\Phi}_{kk}$, are $O_p(n^{2d_k})$. Now,

$$\sum_{j\in N_k} \lambda_j \le \sum_{j\in N_k\cup\cdots\cup N_s} \lambda_j \le \sum_{v=k}^{s} \operatorname{tr}\{\mathbf{\Phi}_{vv}\} = O_p(n^{2d_k}).$$

See, for example, Theorem 14 of Magnus and Neudecker ([12], page 211). We have $\lambda_j = O_p(n^{2d_k})$ for $j \in N_k$. $\square$

Proof of Lemma 3. Following from Lemma 1, $\mathbf{Q}_n^{(k)}$ converges in distribution to a matrix that is positive definite with probability one. Since an eigenvalue of a matrix is a continuous function of the entries of the matrix, we conclude that $\lambda_{j_k^*}(\mathbf{Q}_n^{(k)})$, the smallest eigenvalue of $\mathbf{Q}_n^{(k)}$, converges in distribution to a random variable that has no mass at zero. To prove $n^{-2d_k}\lambda_{j_k^*} \ge c_k \lambda_{j_k^*}(\mathbf{Q}_n^{(k)})$, we construct another, similar, matrix for $I_m(y_t, y_t)$. Let $\mathcal{C}_s = \mathcal{M}(\mathbf{A}_s)$ and $\mathcal{C}_k$, $k = 0, \ldots, s-1$, be the subspaces such that

$$\mathcal{M}^\perp(\mathbf{A}_{k+1}, \ldots, \mathbf{A}_s) = \mathcal{M}^\perp(\mathbf{A}_k, \ldots, \mathbf{A}_s) \oplus \mathcal{C}_k$$

and $\mathcal{C}_k \perp \mathcal{M}^\perp(\mathbf{A}_k, \ldots, \mathbf{A}_s)$. For $k \in \{0, \ldots, s\}$, let $\mathbf{C}_k$ be a $q \times a_k$ matrix with orthonormal columns such that $\mathcal{M}(\mathbf{C}_k) = \mathcal{C}_k$ and $\mathbf{C} = [\mathbf{C}_0 \;\ldots\; \mathbf{C}_s]$. By this construction, $\mathbf{P} = \mathbf{C}'\mathbf{A}$ is a lower triangular block matrix and $\mathbf{P} I_m(z_t, z_t)\mathbf{P}' = \mathbf{C}' I_m(y_t, y_t)\mathbf{C} := \mathbf{W}$ is similar to $I_m(y_t, y_t)$. Let $\mathbf{P}^{(k)}$, $\mathbf{W}^{(k)}$ and $\tilde{\mathbf{d}}_n^{(k)}$, $k = 0, 1, \ldots, s$, be the leading $j_k^* \times j_k^*$ principal submatrices of $\mathbf{P}$, $\mathbf{W}$ and $\mathbf{d}_n$, respectively. Also, let $\tilde{z}_t^{(k)} = (u_t^{(0)}, \ldots, u_t^{(k)})$, $k = 0, 1, \ldots, s$. Note that $\mathbf{W}^{(k)} = \mathbf{P}^{(k)\prime} I_m(\tilde{z}_t^{(k)}, \tilde{z}_t^{(k)})\mathbf{P}^{(k)\prime}$. By Corollary 2.2.1 of Anderson and Das Gupta [1],

$$\lambda_{j_k^*}(\mathbf{W}^{(k)}) \ge \lambda_{j_k^*}\{I_m(\tilde{z}_t^{(k)}, \tilde{z}_t^{(k)})\}\lambda_{j_k^*}(\mathbf{P}^{(k)\prime}\mathbf{P}^{(k)}) = c_k \lambda_{j_k^*}\{I_m(\tilde{z}_t^{(k)}, \tilde{z}_t^{(k)})\}$$

and

$$\begin{aligned}\lambda_{j_k^*}\{I_m(\tilde{z}_t^{(k)}, \tilde{z}_t^{(k)})\} &= \lambda_{j_k^*}\{\tilde{\mathbf{d}}_n^{(k)}\mathbf{Q}_n^{(k)}\tilde{\mathbf{d}}_n^{(k)}\} \ge \lambda_{j_k^*}(\tilde{\mathbf{d}}_n^{(k)})\lambda_{j_k^*}(\mathbf{Q}_n^{(k)})\lambda_{j_k^*}(\tilde{\mathbf{d}}_n^{(k)})\\ &= n^{2d_k}\lambda_{j_k^*}(\mathbf{Q}_n^{(k)}).\end{aligned}$$

Applying the Sturmian separation theorem ([15], page 64), we have

$$\lambda_{j_k^*} = \lambda_{j_k^*}(\mathbf{W}) \ge \lambda_{j_k^*}(\mathbf{W}^{(k)}) \ge c_k \lambda_{j_k^*}\{I_m(\tilde{z}_t^{(k)}, \tilde{z}_t^{(k)})\} \ge c_k n^{2d_k}\lambda_{j_k^*}(\mathbf{Q}_n^{(k)}). \quad \square$$

8.2. *Proofs for Section 5.*

PROOF OF LEMMA 4. Since $\mathbf{X}_k^{\perp}(\mathbf{H}) = [\mathbf{X}_0(\mathbf{H}) \ \ldots \ \mathbf{X}_{k-1}(\mathbf{H}) \ \mathbf{X}_{k+1}(\mathbf{H}) \ \ldots \ \mathbf{X}_q(\mathbf{H})]$, we have

$$\begin{aligned}\|\sin\Theta\{\mathcal{M}(\mathbf{X}_k(\mathbf{H})),\mathcal{M}(\mathbf{X}_k)\}\|_F &\le \|(\mathbf{X}_k^{\perp}(\mathbf{H}))^*\mathbf{X}_k\|_F \le \sum_{\ell=0,\ell\neq k}^{s}\|(\mathbf{X}_\ell(\mathbf{H}))^*\mathbf{X}_k\|_F \\ &= O_p\Big(\max_{\ell\neq k} n^{-|d_k-d_\ell|}\Big) = O_p(n^{-\alpha_k}),\end{aligned}$$

by Lemma 7. □

PROOF OF LEMMA 5. For $k = 1,\ldots,s-1$, we have

$$\begin{aligned}&P(\mathcal{M}(\mathbf{X}_k(\mathbf{H})) = \mathcal{B}_k) \\ &\quad = P\bigg(\bigg\{\mathcal{M}\mathbf{X}_k(\mathbf{H}) \cap \bigoplus_{\ell\le k-1}\mathcal{B}_\ell = \mathbf{0}\bigg\} \cap \bigg\{\mathcal{M}\mathbf{X}_k(\mathbf{H}) \cap \bigoplus_{\ell\ge k+1}\mathcal{B}_\ell = \mathbf{0}\bigg\}\bigg).\end{aligned}$$

Hence,

$$\begin{aligned}&P(\mathcal{M}(\mathbf{X}_k(\mathbf{H})) \neq \mathcal{B}_k) \\ &\quad = P\bigg(\bigg\{\mathcal{M}\mathbf{X}_k(\mathbf{H}) \cap \bigoplus_{\ell\le k-1}\mathcal{B}_\ell \neq \mathbf{0}\bigg\} \cup \bigg\{\mathcal{M}\mathbf{X}_k(\mathbf{H}) \cap \bigoplus_{\ell\ge k+1}\mathcal{B}_\ell \neq \mathbf{0}\bigg\}\bigg) \\ &\quad \le P\bigg(\mathcal{M}\mathbf{X}_k(\mathbf{H}) \cap \bigoplus_{\ell\le k-1}\mathcal{B}_\ell \neq \mathbf{0}\bigg) + P\bigg(\mathcal{M}\mathbf{X}_k(\mathbf{H}) \cap \bigoplus_{\ell\ge k+1}\mathcal{B}_\ell \neq \mathbf{0}\bigg) \\ &\quad = O(n^{-2d_{k-1}+2d_k} + n^{-2d_k+2d_{k+1}}),\end{aligned}$$

by Lemma 10. Similarly,

$$P(\mathcal{M}(\mathbf{X}_0(\mathbf{H})) \neq \mathcal{B}_0) = O_p(n^{-2d_0+2d_1})$$

and

$$P(\mathcal{M}(\mathbf{X}_s(\mathbf{H})) \neq \mathcal{B}_s) = O_p(n^{-2d_{s-1}+2d_s}).$$

We have thus completed the proof. □

We will need the following lemma for the proof of Lemma 7. First we write $\mathbf{U}$ and $\mathbf{V}$ defined in Lemma 1 as

$$\mathbf{U} = [\mathbf{U}_0' \ \ldots \ \mathbf{U}_s']', \qquad \mathbf{V} = [\mathbf{V}_0' \ \ldots \ \mathbf{V}_s']',$$

where $\mathbf{U}_k$ and $\mathbf{V}_k$ are $a_k \times m$ matrices.

LEMMA 6. *Let* $\mathbf{K} = \mathrm{diag}(\mathbf{B}_0'\mathbf{A}_0, \ldots, \mathbf{B}_s'\mathbf{A}_s)$. *Then*

$$\mathbf{d}_n^{-1}\mathbf{\Phi}\mathbf{d}_n^{-1} \xrightarrow{D} \mathbf{K}(\mathbf{U}\mathbf{U}' + \mathbf{V}\mathbf{V}')\mathbf{K}',$$

*where* $\mathbf{d}_n$ *is defined as in Lemma* 1. *Furthermore,*

$$\mathbf{d}_n^{-1}\mathbf{\Phi}_D\mathbf{d}_n^{-1} \xrightarrow{D} \mathbf{K}(\mathbf{U}\mathbf{U}' + \mathbf{V}\mathbf{V}')_D\mathbf{K}',$$

*where* $(\mathbf{U}\mathbf{U}' + \mathbf{V}\mathbf{V}')_D = \mathrm{diag}(\mathbf{U}_0\mathbf{U}_0' + \mathbf{V}_0\mathbf{V}_0', \ldots, \mathbf{U}_s\mathbf{U}_s' + \mathbf{V}_s\mathbf{V}_s')$ *and* $\mathbf{K}(\mathbf{U}\mathbf{U}' + \mathbf{V}\mathbf{V}')_D\mathbf{K}'$ *is positive definite and has distinctive eigenvalues with probability* 1.

PROOF. We write $\mathbf{\Phi} = \mathbf{K}I_m(z_t, z_t)\mathbf{K} + \mathbf{R}$, where $\mathbf{R}$ is a symmetric matrix with its $(k, \ell)$th entry given by

$$\begin{aligned}\mathbf{R}_{k\ell} = {} & \mathbf{B}_k'\mathbf{A}_k I_m(u_t^{(k)}, z_t^{(\ell+1)})\mathbf{A}^{(\ell+1)\prime}\mathbf{B}_\ell + \mathbf{B}_k'\mathbf{A}^{(k+1)} I_m(z_t^{(k+1)}, u_t^{(\ell)})\mathbf{A}_\ell'\mathbf{B}_\ell \\ & + \mathbf{B}_k'\mathbf{A}^{(k+1)} I_m(z_t^{(k+1)}, z_t^{(\ell+1)})\mathbf{A}^{(\ell+1)\prime}\mathbf{B}_\ell\end{aligned}$$

for $k \le \ell$, $\ell = 0, 1, \ldots, (s-1)$, $\mathbf{R}_{ks} = \mathbf{B}_k'\mathbf{A}^{(k+1)} I_m(z_t^{(k+1)}, z_t^{(s)})\mathbf{A}_s'\mathbf{B}_s$ for $k < s$ and $\mathbf{R}_{ss} = 0$. Thus,

$$\mathbf{d}_n^{-1}\mathbf{\Phi}\mathbf{d}_n^{-1} = \mathbf{d}_n^{-1}\mathbf{K}I_m(z_t, z_t)\mathbf{K}'\mathbf{d}_n^{-1} + \mathbf{d}_n^{-1}\mathbf{R}\mathbf{d}_n^{-1}.$$

By Lemma 1,

$$\mathbf{d}_n^{-1}\mathbf{K}I_m(z_t, z_t)\mathbf{K}'\mathbf{d}_n^{-1} \xrightarrow{D} \mathbf{K}(\mathbf{U}\mathbf{U}' + \mathbf{V}\mathbf{V}')\mathbf{K}'$$

and the $(k, \ell)$th entry of $\mathbf{d}_n^{-1}\mathbf{R}\mathbf{d}_n^{-1}$ is

$$n^{-d_k - d_\ell}\mathbf{R}_{k\ell} = O_p(n^{d_{\ell+1} - d_\ell} + n^{d_{k+1} - d_k}) = o_p(1) \qquad \text{for } k \le \ell.$$

We have proved the first limiting distribution of the lemma. It follows that the $k$th diagonal block of $\mathbf{\Phi}$ has the limiting distribution

$$n^{-2d_k}\mathbf{\Phi}_{kk} \xrightarrow{D} \mathbf{K}_{kk}(\mathbf{U}_k\mathbf{U}_k' + \mathbf{V}_k\mathbf{V}_k')\mathbf{K}_{kk}'$$

and $\mathring{\mathbf{\Phi}}_{kk}$ is positive definite, having distinctive eigenvalues with probability 1 by Okamoto [14]. □

LEMMA 7. $\|\mathbf{X}_\ell^*(\mathbf{H})\mathbf{X}_k\|_F = O_p(n^{-|d_k - d_\ell|})$ *for all* $\ell, k \in \{0, 1, \ldots, s\}$ *with* $\ell \ne k$.

PROOF. Since $\|\mathbf{X}_\ell^*(\mathbf{H})\mathbf{X}_k\|_F = \|\mathbf{X}_\ell^*(\mathbf{\Phi}_D)\mathbf{B}'\mathbf{B}\mathbf{X}_k(\mathbf{\Phi})\|_F = \|\mathbf{X}_\ell^*(\mathbf{\Phi}_D)\mathbf{X}_k(\mathbf{\Phi})\|_F$, we prove this lemma by showing that

$$\|\mathbf{X}_\ell^*(\mathbf{\Phi}_D)\mathbf{X}_k(\mathbf{\Phi})\|_F = O_p(n^{-|d_k - d_\ell|}).$$

Let $\mathbf{\Lambda} = \mathrm{diag}\{\lambda_j, j = 1, \ldots, q\}$ and $\Lambda^{(k)} = \{\lambda_j, j \in N_k\}$. We define $\mathbf{\Lambda}(\mathbf{\Phi}_D)$ and $\mathbf{\Lambda}^{(k)}(\mathbf{\Phi}_D)$ similarly for $\mathbf{\Phi}_D$. We will use the bound for the error in two subspaces within the nonzero space from Theorem 4.1 of Barlow and Slapničar [2] (which can be shown to apply in our context with probability one), that is,

$$\|\mathbf{X}_\ell^*(\mathbf{\Phi}_D)\mathbf{X}_k(\mathbf{\Phi})\|_F \le \frac{\|\mathbf{\Lambda}^{-1/2}(\mathbf{\Phi}_D)\mathbf{X}^*(\mathbf{\Phi}_D)\Delta\mathbf{\Phi}\mathbf{X}(\mathbf{\Phi})\mathbf{\Lambda}^{-1/2}\|_F}{\mathrm{relgap}(\Lambda^{(\ell)}(\mathbf{\Phi}_D),\Lambda^{(k)})},$$

where

$$\mathrm{relgap}(\Lambda^{(\ell)}(\mathbf{\Phi}_D),\Lambda^{(k)}) = \min_{i\in N_k, j\in N_\ell}\left|\frac{\lambda_i(\mathbf{\Phi})-\lambda_j(\mathbf{\Phi}_D)}{\lambda_i^{1/2}(\mathbf{\Phi}_D)\lambda_j^{1/2}(\mathbf{\Phi})}\right|.$$

It is sufficient to show that

$$\|\mathbf{\Lambda}^{-1/2}(\mathbf{\Phi}_D)\mathbf{X}^*(\mathbf{\Phi}_D)\Delta\mathbf{\Phi}\mathbf{X}(\mathbf{\Phi})\mathbf{\Lambda}^{-1/2}\|_F = O_p(1) \tag{16}$$

and

$$\frac{1}{\mathrm{relgap}(\Lambda^{(\ell)}(\mathbf{\Phi}_D),\Lambda^{(k)})} = O_p(n^{-|d_k-d_\ell|}). \tag{17}$$

By Lemmas 2, 3 and 6, $\mathrm{relgap}(\Lambda^{(\ell)}(\mathbf{\Phi}_D),\Lambda^{(k)}) = O_p(n^{|d_k-d_\ell|})$ and $n^{-|d_k-d_\ell|} \times \mathrm{relgap}(\Lambda^{(\ell)}(\mathbf{\Phi}_D),\Lambda^{(k)}) \ge \varsigma_{\ell,k}$, where $\varsigma_{\ell,k}$ is a random variable that has no mass at 0. We thus have (17). We next prove (16). Note that by Lemmas 1 and 6,

$$\mathbf{d}_n\mathbf{\Phi}^{-1}\mathbf{d}_n \xrightarrow{D} \mathbf{K}'^{-1}(\mathbf{U}\mathbf{U}'+\mathbf{V}\mathbf{V}')^{-1}\mathbf{K}^{-1}.$$

Hence, $\mathbf{d}_n\mathbf{X}(\mathbf{\Phi})\mathbf{\Lambda}^{-1/2} = O_p(1)$ since $\mathbf{d}_n\mathbf{\Phi}^{-1}\mathbf{d}_n = \mathbf{d}_n\mathbf{X}(\mathbf{\Phi})\mathbf{\Lambda}^{-1/2}\mathbf{\Lambda}^{-1/2}\mathbf{X}'(\mathbf{\Phi})\mathbf{d}_n = O_p(1)$. Similarly, $\mathbf{\Lambda}^{-1/2}(\mathbf{\Phi}_D)\mathbf{X}^*(\mathbf{\Phi}_D)\mathbf{d}_n = O_p(1)$. We have

$$\begin{aligned}
&\|\mathbf{\Lambda}^{-1/2}(\mathbf{\Phi}_D)\mathbf{X}^*(\mathbf{\Phi}_D)\Delta\mathbf{\Phi}\mathbf{X}(\mathbf{\Phi})\mathbf{\Lambda}^{-1/2}\|_F \\
&\quad = \|\mathbf{\Lambda}^{-1/2}(\mathbf{\Phi}_D)\mathbf{X}^*(\mathbf{\Phi}_D)\mathbf{d}_n\mathbf{d}_n^{-1}\Delta\mathbf{\Phi}\mathbf{d}_n^{-1}\mathbf{d}_n\mathbf{X}(\mathbf{\Phi})\mathbf{\Lambda}^{-1/2}\|_F \\
&\quad \le \|\mathbf{\Lambda}^{-1/2}(\mathbf{\Phi}_D)\mathbf{X}^*(\mathbf{\Phi}_D)\mathbf{d}_n\|_F\|\mathbf{d}_n^{-1}\Delta\mathbf{\Phi}\mathbf{d}_n^{-1}\|_F\|\mathbf{d}_n\mathbf{X}(\mathbf{\Phi})\mathbf{\Lambda}^{-1/2}\|_F \\
&\quad = O_p(1),
\end{aligned}$$

by Lemma 6. Hence, $\|\mathbf{X}_\ell^*(\mathbf{\Phi}_D)\mathbf{X}_k(\mathbf{\Phi})\|_F = O_p(n^{-|d_k-d_\ell|})$. □

We need the following two lemmas for the proof of Lemma 9:

Lemma 8. *Under Assumption 1, there exists a finite constant $C$ not depending on $n$ such that for all sufficiently large $n$,*

$$\mathbb{E}[\lambda_1^2(\mathbf{Q}_n^{-1})] < C.$$

PROOF. Note that

$$\mathbf{Q}_n = (\mathbf{U}_n, \mathbf{V}_n)(\mathbf{U}_n, \mathbf{V}_n)',$$

where $\mathbf{U}_n$ and $\mathbf{V}_n$ are defined in Eq. (10). Let

$$\mathcal{T}(W_n) = \lambda_1^2(\mathbf{Q}_n^{-1}),$$

where $W_n=$ vec$(\mathbf{U}_n, \mathbf{V}_n)$. By Assumption 1, $W_n \sim N(0, \boldsymbol{\Xi}_n)$, where $\boldsymbol{\Xi}_n = cov(W_n)$ and $\boldsymbol{\Xi}_n \to \boldsymbol{\Xi}$, the covariance matrix of vec$(\mathbf{U}, \mathbf{V})$ in Lemma 1. It was shown in [4] that $\boldsymbol{\Xi}$ is positive definite. Thus, for all sufficiently large $n$, $\boldsymbol{\Xi}_n$ is invertible and $\lambda_1(\boldsymbol{\Xi}_n) \to \lambda_1(\boldsymbol{\Xi}) > 0$.

For all sufficiently large $n$,

$$\mathbb{E}_{\boldsymbol{\Xi}_n}[\mathcal{T}(W_n)] = (2\pi)^{-mq}|\boldsymbol{\Xi}_n|^{-1/2} \int_{\mathbb{R}^{2mq}} \mathcal{T}(x) e^{-x'\boldsymbol{\Xi}_n^{-1}x/2}\,dx.$$

Since $x'\boldsymbol{\Xi}_n^{-1}x' \geq x'x/\lambda_1(\boldsymbol{\Xi}_n)$, we have

$$e^{-x'\boldsymbol{\Xi}_n^{-1}x/2} \leq e^{-x'x/2\lambda_1(\boldsymbol{\Xi}_n)}.$$

Since $\lambda_1(\boldsymbol{\Xi}_n) \to \lambda_1(\boldsymbol{\Xi}) > 0$ and since $|\boldsymbol{\Xi}_n|^{-1/2} \to |\boldsymbol{\Xi}|^{-1/2} > 0$, there exist constants $C_1 > 0$ and $C_2 > 0$ such that for all sufficiently large $n$,

$$\mathbb{E}_{\boldsymbol{\Xi}_n}[\mathcal{T}(W_n)] \leq C_1 \int_{\mathbb{R}^{2mq}} \mathcal{T}(x) e^{-C_2 x'x/2}\,dx = C,$$

a finite constant which does not depend on $n$. The above integral is the second moment of the largest eigenvalue of an inverse Wishart matrix and is hence bounded by a finite constant [21], in view of our assumption that $m > q + 3$. $\square$

LEMMA 9. *Define $E_{k\ell}$ to be an event, $E_{k\ell} = \{\lambda_{a_k}(\mathbf{\Phi}_{kk}) > \lambda_1(\mathbf{\Phi}_{\ell\ell})\}$, $0 \leq k < \ell \leq s$. Then under Assumption 1,*

$$P(E_{k\ell}^c) = O(n^{-2d_k+2d_\ell}).$$

PROOF. For $\ell > k$, $\ell = 1, \dots, s$, we have, by Chebyshev's inequality and the Cauchy–Schwarz inequality,

$$\begin{aligned} P(E_{k\ell}^c) &= P\{\lambda_{a_k}(\mathbf{\Phi}_{kk}) \leq \lambda_1(\mathbf{\Phi}_{\ell\ell})\} \\ &= P\left\{n^{-2d_k+2d_\ell} \cdot \frac{n^{-2d_\ell}\lambda_1(\mathbf{\Phi}_{\ell\ell})}{n^{-d_k}\lambda_{a_k}(\mathbf{\Phi}_{kk})} \geq 1\right\} \\ &\leq n^{-2d_k+2d_\ell}\mathbb{E}^{1/2}[\lambda_1^2(n^{-2d_\ell}\mathbf{\Phi}_{\ell\ell})]\mathbb{E}^{1/2}[\lambda_1^2(n^{2d_k}\mathbf{\Phi}_{kk}^{-1})]. \end{aligned}$$

Since $\mathbb{E}[\lambda_1^2(n^{-2d_\ell}\mathbf{\Phi}_{\ell\ell})] \leq \mathbb{E}[\mathrm{trace}^2(n^{-2d_\ell}\mathbf{\Phi}_{\ell\ell})] < C$ by Assumption 1 and Lemma 1, the lemma follows if

$$\mathbb{E}[\lambda_1^2(n^{2d_k}\mathbf{\Phi}_{kk}^{-1})] < C$$

for all sufficiently large $n$. Let $\mathbf{J}=\mathbf{B}_k'\mathbf{A}^{(k)}\mathbf{d}_n^{(k)}$, $\mathbf{d}_n^{(k)}=\operatorname{diag}(n^{d_k},\dots,n^{d_k},\dots,n^{d_{u_s}},\dots,n^{d_{u_s}})'$. We write

$$\boldsymbol{\Phi}_{kk}=\mathbf{J}\{(\mathbf{d}_n^{(k)})^{-1}I_m(z_t^{(k)},z_t^{(k)})(\mathbf{d}_n^{(k)})^{-1}\}\mathbf{J}'.$$

We will use the inequality of Exercise 19 on page 238 of Magnus and Neudecker [12]. That is,

$$\boldsymbol{\Phi}_{kk}^{-1}\le(\mathbf{JJ}')^{-1}\mathbf{J}\{\mathbf{d}_n^{(k)}[I_m(z_t^{(k)},z_t^{(k)})]^{-1}\mathbf{d}_n^{(k)}\}\mathbf{J}'(\mathbf{JJ}')^{-1}.$$

It follows that

$$\begin{aligned}\operatorname{trace}\boldsymbol{\Phi}_{kk}^{-1}&\le\lambda_1\{\mathbf{d}_n^{(k)}I_m^{-1}(z_t^{(k)},z_t^{(k)})\mathbf{d}_n^{(k)}\}\operatorname{trace}\{(\mathbf{JJ}')^{-1}\mathbf{JJ}'(\mathbf{JJ}')^{-1}\}\\&=\lambda_1\{\mathbf{d}_n^{(k)}I_m^{-1}(z_t^{(k)},z_t^{(k)})\mathbf{d}_n^{(k)}\}\operatorname{trace}\{(\mathbf{JJ}')^{-1}\}.\end{aligned}$$

Since there exists a finite constant $C$ such that $E[\lambda_1^2\{\mathbf{d}_n^{(k)}I_m^{-1}(z_t^{(k)},z_t^{(k)})\mathbf{d}_n^{(k)}\}]<C$ for all sufficiently large $n$ by Lemma 8, we complete the proof by showing that

$$\operatorname{trace}\{(\mathbf{JJ}')^{-1}\}=O(n^{-2d_k}).$$

We write

$$\mathbf{JJ}'=n^{2d_k}\mathbf{B}_k'\mathbf{A}_k\mathbf{A}_k'\mathbf{B}_k+\mathbf{B}_k'\mathbf{A}^{(k+1)}\mathbf{d}_n^{(k+1)}\mathbf{d}_n^{(k+1)}\mathbf{A}^{(k+1)}\mathbf{B}_k\mathbf{1}_{\{k<s\}}.$$

Since both matrices on the right-hand side are symmetric and positive definite, $\lambda_{a_k}\mathbf{JJ}'\ge\lambda_{a_k}[n^{2d_k}\mathbf{B}_k'\mathbf{A}_k\mathbf{A}_k'\mathbf{B}_k]$, and we have

$$\lambda_1\{(\mathbf{JJ}')^{-1}\}\le n^{-2d_k}\lambda_1\{[\mathbf{B}_k'\mathbf{A}_k\mathbf{A}_k'\mathbf{B}_k]^{-1}\}=O(n^{-2d_k}).\qquad\square$$

LEMMA 10. *Under Assumption 1,*

$$P\Big\{\mathcal{M}\mathbf{X}_k(\mathbf{H})\cap\bigoplus_{\ell\le h_1}\mathcal{B}_j\ne\mathbf{0}\Big\}=O(n^{-2d_{h_1}+2d_k})\tag{18}$$

*for* $h_1<k$, $k=1,\dots,s$ *and*

$$P\Big\{\mathcal{M}\mathbf{X}_k(\mathbf{H})\cap\bigoplus_{\ell\ge h_2}\mathcal{B}_j\ne\mathbf{0}\Big\}=O(n^{-2d_k+2d_{h_2}})\tag{19}$$

*for* $h_2>k$, $k=0,\dots,s-1$.

PROOF. Since $\mathbf{H}=\mathbf{B}\boldsymbol{\Phi}_D\mathbf{B}'$, we have $\mathbf{X}_\ell(\mathbf{H})=\mathbf{B}\mathbf{X}_\ell(\boldsymbol{\Phi}_D)$. Since $\boldsymbol{\Phi}_D$ is a block diagonal matrix,

$$\lambda_i(\boldsymbol{\Phi}_D)\in\{\lambda_j(\boldsymbol{\Phi}_{kk})|k=0,\dots,s,\ j=1,\dots,a_k\}$$

and for $\lambda_i(\boldsymbol{\Phi}_D)$ such that $\lambda_i(\boldsymbol{\Phi}_D)=\lambda_j(\boldsymbol{\Phi}_{kk})$,

$$\chi_i(\boldsymbol{\Phi}_D)=(0,\dots,0,\chi_j'(\boldsymbol{\Phi}_{kk}),0,\dots0)',$$

that is, the first $j^*_{k-1}$ entries are all zero. Define $E_{h\ell}$ to be an event, $E_{h\ell} = \{\lambda_{a_h}(\mathbf{\Phi}_{hh}) > \lambda_1(\mathbf{\Phi}_{\ell\ell})\}$, $0 \le h < \ell \le s$. We first prove (18).

$$P\Big\{\mathcal{M}\mathbf{X}_k(\mathbf{H}) \cap \bigoplus_{\ell \le h_1} \mathcal{B}_\ell \ne \mathbf{0}\Big\} = P(\mathbf{X}_k(\mathbf{\Phi}_D) \ne [\mathbf{0}\ \mathbf{Y}]'),$$

where the 0 in $[\mathbf{0}\ \mathbf{Y}]'$ has dimension $j^*_{h_1} \times a_k$ and $\mathbf{Y}$ has full rank. We have for $h_1 < k$, $k = 1, \dots, s$,

$$\begin{aligned} P(\mathbf{X}_k(\mathbf{\Phi}_D) \ne [\mathbf{0}\ \mathbf{Y}]') = P\Big(\bigcup_{\ell:\ell\le h_1} E^c_{\ell k}\Big) \le \sum_{\ell:\ell \le h_1} P(E^c_{\ell k}) = O\Big(\sum_{\ell:\ell\le h_1} n^{-2d_\ell + 2d_k}\Big) \\ = O(n^{-2d_{h_1}+2d_k}), \end{aligned}$$

by Lemma 9. Similarly, for (19),

$$P\Big\{\mathcal{M}\mathbf{X}_k(\mathbf{H}) \cap \oplus_{\ell \ge h_2} \mathcal{B}_\ell \ne \mathbf{0}\Big\} = P(\mathbf{X}_k(\mathbf{\Phi}_D) \ne [\mathbf{Z}\ \mathbf{0}]'),$$

where the 0 in $[\mathbf{Z}\ \mathbf{0}]'$ has dimension $(q - j^*_{h_2}) \times a_k$, and $\mathbf{Z}$ has full rank. We have for $h_2 > k$, $k = 0, \dots, s-1$,

$$\begin{aligned} P(\mathbf{X}_k(\mathbf{\Phi}_D) \ne [\mathbf{Z}\ \mathbf{0}]') = P\Big(\bigcup_{\ell:\ell\ge h_2} E^c_{k\ell}\Big) \le \sum_{\ell:\ell \ge h_2} P(E^c_{k\ell}) = O\Big(\sum_{\ell:\ell\ge h_2} n^{-2d_k + 2d_\ell}\Big) \\ = O(n^{-2d_k+2d_{h_2}}). \quad \square \end{aligned}$$

8.3. *Proofs for Sections 6 and 7.* In this section, we will use the following decomposition and notation for the proofs. We write

$$b'\mathbf{A}I_{zz}(\omega_j)\mathbf{A}'b - b'\mathbf{A}\mathbf{f}(\omega_{\tilde{j}})\mathbf{A}'b = b'\mathbf{A}\mathbf{R}(\omega_j)\mathbf{A}'b + b'\mathbf{A}\mathbf{S}(\omega_j)\mathbf{A}'b, \tag{20}$$

where

$$\mathbf{R}(\omega_j) = I_{zz}(\omega_j) - \mathbf{\Psi}(\omega_{\tilde{j}})I_{\varepsilon\varepsilon}(\omega_j)\mathbf{\Psi}^*(\omega_{\tilde{j}})$$

and

$$\mathbf{S}(\omega_j) = \mathbf{\Psi}(\omega_{\tilde{j}})I_{\varepsilon\varepsilon}(\omega_j)\mathbf{\Psi}^*(\omega_{\tilde{j}}) - \mathbf{f}(\omega_{\tilde{j}}).$$

We will also use the following notation:

$$\begin{aligned} \mathcal{L}_{m_n}(d) &= \frac{1}{m_n}\sum_{j=1}^{m_n} \omega_{\tilde{j}}^{2d} b'\mathbf{A}\mathbf{R}(\omega_j)\mathbf{A}'b, \\ \mathcal{M}_{m_n}(d) &= \frac{1}{m_n}\sum_{j=1}^{m_n} \omega_{\tilde{j}}^{2d} b'\mathbf{A}\mathbf{S}(\omega_j)\mathbf{A}'b, \\ \mathcal{F}_{m_n}(d) &= \frac{1}{m_n}\sum_{j=1}^{m_n} \omega_{\tilde{j}}^{2d} b'\mathbf{A}\mathbf{f}(\omega_j)\mathbf{A}'b. \end{aligned} \tag{21}$$

8.3.1. *Proof of Theorem 2*. Here and in the following subsections, suppose that $a \in N_k$, where $k \in \{0, \ldots, s\}$. Write $\hat{d}_{aa} = \hat{d}_k$. Note that $d_{aa} = d_k$. For $1/4 > \delta > 0$, let $\mathcal{N}_\delta = \{d : |d - d_k| < \delta\}$. Then for $S(d) = R(d) - R(d_k)$, we have

$$P(|\hat{d}_k - d_k| \geq \delta) = P(\hat{d}_k \in \mathcal{N}_\delta{}^c \cap \Theta)$$
$$= P\Big(\inf_{\mathcal{N}_\delta{}^c \cap \Theta} R(d) \leq \inf_{\mathcal{N}_\delta \cap \Theta} R(d)\Big) \leq P\Big(\inf_{\mathcal{N}_\delta{}^c \cap \Theta} S(d) \leq 0\Big).$$

Define $\Theta_1 = \{d : \Delta \leq d \leq \Delta_2\}$, where $\Delta = \Delta_1$ when $d_k < 1/2 + \Delta_1$ and $d_k \geq \Delta > d_k - 1/2$ otherwise. Note that $d - d_k > -1/2$ for all $d \in \Theta_1$. When $d_k \geq 1/2 + \Delta_1$, define $\Theta_2 = \{d : \Delta_1 \leq d < \Delta\}$ and otherwise take $\Theta_2$ to be empty. Hence,

$$P(|\hat{d}_k - d_k| \geq \delta) \leq P\Big(\inf_{\mathcal{N}_\delta{}^c \cap \Theta_1} S(d) \leq 0\Big) + P\Big(\inf_{\Theta_2} S(d) \leq 0\Big) = o(1),$$

by Lemmas 11 and 12 below. □

LEMMA 11. *Under the assumptions of Theorem 2*, $P(\inf_{\mathcal{N}_\delta^c \cap \Theta_1} S(d) \leq 0) = o(1)$.

PROOF. Let

$$U(d) = 2(d - d_k) - \log\{2(d - d_k) + 1\}$$

and

$$T(d) = \log \frac{\widehat{G}(d_k)}{\mathcal{G}} - \log \frac{\widehat{G}(d)}{G(d)} - \log\left\{\frac{2(d - d_k) + 1}{m_n} \sum_{j=1}^{m_n} \left(\frac{\tilde{j}}{m_n}\right)^{2(d-d_k)}\right\}$$
$$+ 2(d - d_k)\left\{\frac{1}{m_n} \sum_{j=1}^{m_n} \log \tilde{j} - (\log m_n - 1)\right\},$$

where $\mathcal{G} = b'\mathbf{A}_k \mathbf{f}^\dagger(0)\mathbf{A}_k' b$, as in (31), and

$$G(d) = \mathcal{G} \frac{1}{m_n} \sum_{j=1}^{m_n} \omega_{\tilde{j}}^{2(d-d_k)}.$$

Then $S(d) = U(d) - T(d)$. We have

$$P\Big(\inf_{\mathcal{N}_\delta^c \cap \Theta_1} S(d) \leq 0\Big) \leq P\Big(\inf_{\mathcal{N}_\delta^c \cap \Theta_1} U(d) \leq \sup_{\Theta_1} |T(d)|\Big).$$

Following the same arguments as those on page 1635 of Robinson [17], it is sufficient to show that

$$\sup_{\Theta_1} \left|\frac{\widehat{G}(d) - G(d)}{G(d)}\right| = o_p(1).$$

Note that by Lemma 20, we have

$$G(d) = C\mathcal{G}\omega_{m_n}^{2(d-d_k)} \geq C(1-\varepsilon_k)\omega_{m_n}^{2(d-d_k)}, \tag{22}$$

where $\varepsilon_k = O_p(n^{-\alpha_k})$. By Lemma 21, for $d \in \Theta_1$,

$$\begin{aligned} |\widehat{G}(d) - G(d)| &= \mathcal{L}_{m_n}(d) + \mathcal{M}_{m_n}(d) + \mathcal{F}_{m_n}(d) - \mathcal{G}\frac{1}{m_n}\sum_{j=1}^{m_n}\omega_{\tilde{j}}^{2(d-d_k)} \\ &= o_p(\omega_{m_n}^{2d-2d_k} m_n^{-\varepsilon}), \end{aligned} \tag{23}$$

where $\mathcal{L}_{m_n}$, $\mathcal{M}_{m_n}$ and $\mathcal{F}_{m_n}$ are defined in (21). We have thus completed the proof. □

Lemma 12. *Under the assumptions of Theorem 2,* $P(\inf_{\Theta_2} S(d) \leq 0) = o(1)$.

Proof. Following from the proof on pages 1638–1639 of Robinson [17], we write

$$S(d) = \log\{\widehat{D}(d)/\widehat{D}(d_k)\},$$

where

$$\widehat{D}(d) = \frac{1}{m_n}\sum_{j=1}^{m_n}\left(\frac{\tilde{j}}{e^\nu}\right)^{2(d-d_k)} \tilde{j}^{-2d_k} I_{vv}(\omega_j) \quad \text{and} \quad \nu = \frac{1}{m_n}\sum_{j=1}^{m_n}\log\tilde{j}.$$

Note that $e^\nu \sim m_n/e$. Denote

$$\alpha_j = \begin{cases} \left(\dfrac{\tilde{j}}{e^\nu}\right)^{2(\Delta-d_k)} \sim \left(\dfrac{ej}{m_n}\right)^{2(\Delta-d_k)}, & 1 \leq j \leq e^\nu, \\ \left(\dfrac{\tilde{j}}{e^\nu}\right)^{2(\Delta_1-d_k)} \sim \left(\dfrac{ej}{m_n}\right)^{2(\Delta_1-d_k)}, & e^\nu < j < m_n. \end{cases} \tag{24}$$

By choosing $\Delta < d_k - \frac{1}{2} + \frac{1}{4e}$ so that $m_n^{-1}\sum_{j=1}^{m_n}(\alpha_j - 1) \geq 1$ for all sufficiently large $m_n$, we have

$$\begin{aligned} P\Big(\inf_{\Theta_2} S(d) \leq 0\Big) &\leq P\left(\frac{1}{m_n}\sum_{j=1}^{m_n}(\alpha_j - 1)\tilde{j}^{-2d_k} I_{vv}(\omega_j) \leq 0\right) \\ &= P\left(\frac{1}{m_n}\sum_{j=1}^{m_n}(\alpha_j - 1)\frac{I_{vv}(\omega_j)}{\mathcal{G}\omega_{\tilde{j}}^{-2d_k}} \leq 0\right) \\ &\leq P\left(\left|\frac{1}{m_n}\sum_{j=1}^{m_n}(\alpha_j - 1)\left(\frac{I_{vv}(\omega_j)}{\mathcal{G}\omega_{\tilde{j}}^{-2d_k}} - 1\right)\right| \geq 1\right). \end{aligned}$$

Now, by (20),

$$\frac{1}{m_n}\sum_{j=1}^{m_n}(\alpha_j-1)\left(\frac{I_{vv}(\omega_j)}{\mathcal{G}\omega_{\tilde{j}}^{-2d_k}}-1\right)=\frac{1}{m_n}\sum_{j=1}^{m_n}(\alpha_j-1)\left(\frac{I_{vv}(\omega_j)}{\mathcal{G}\omega_{\tilde{j}}^{-2d_k}}-\frac{I_{vv}(\omega_j)}{b'\mathbf{A}\mathbf{f}(\omega_{\tilde{j}})\mathbf{A}'b}\right)$$
$$+\frac{1}{m_n}\sum_{j=1}^{m_n}(\alpha_j-1)\frac{b'\mathbf{A}\mathbf{R}(\omega_j)\mathbf{A}'b}{b'\mathbf{A}\mathbf{f}(\omega_{\tilde{j}})\mathbf{A}'b} \tag{25}$$
$$+\frac{1}{m_n}\sum_{j=1}^{m_n}(\alpha_j-1)\frac{b'\mathbf{A}\mathbf{S}(\omega_j)\mathbf{A}'b}{b'\mathbf{A}\mathbf{f}(\omega_{\tilde{j}})\mathbf{A}'b}.$$

We will show that all three terms in (25) are $o_p(1)$. For the first term, we begin by showing that

$$I_{vv}(\omega_j)=b'\mathbf{A}\mathbf{R}(\omega_j)\mathbf{A}'b+b'\mathbf{A}\mathbf{S}(\omega_j)\mathbf{A}'b+b'\mathbf{A}\mathbf{f}(\omega_j)\mathbf{A}'b=O_p(\omega_j^{-2d_k}). \tag{26}$$

Let $\mathbf{R}_{h\ell}(\omega_j)$ denote the $(h,\ell)$th block of $\mathbf{R}(\omega_j)$. By Lemmas 16 and 18,

$$b'\mathbf{A}_h\mathbf{R}_{h\ell}(\omega_j)\mathbf{A}_\ell'b=\begin{cases}O_p(n^{2d_k-d_h-d_\ell}\omega_j^{-d_h-d_\ell}j^{-\rho/2}), & h,\ell<k,\\ O_p(\omega_j^{-d_h-d_\ell}j^{-\rho/2}), & h,\ell\geq k,\\ O_p(n^{d_k-d_h}\omega_j^{-d_h-d_\ell}j^{-\rho/2}), & h<k,\ell\geq k,\end{cases} \tag{27}$$
$$=O_p(\omega_j^{-2d_k}j^{-\rho/2}).$$

Also, by Lemmas 16 and 19,

$$\begin{aligned}&b'\mathbf{A}\mathbf{S}(\omega_j)\mathbf{A}'b\\ &\quad=O_p\left(\sum_{h,\ell=0}^{k-1}\omega_j^{-d_h-d_\ell}n^{d_{2k}-d_h-d_\ell}+\sum_{h,\ell=k}^{s}\omega_j^{-d_h-d_\ell}+\sum_{h<k,\ell\geq k}\omega_j^{-d_h-d_\ell}n^{d_k-d_h}\right)\\ &\quad=O_p(\omega_j^{-2d_k}(j^{2d_k-2d_{k-1}}+1+j^{d_k-d_{k-1}})).\end{aligned}$$

By (54) in the proof of Lemma 20, $b'\mathbf{A}\mathbf{f}(\omega_j)\mathbf{A}'b=O_p(\omega_j^{-2d_k})$. Thus, the bound in (26) follows. Together with Lemma 20, we have

$$\frac{I_{vv}(\omega_j)}{\mathcal{G}\omega_{\tilde{j}}^{-2d_k}}=O_p(1)$$

and

$$1-\frac{\mathcal{G}\omega_{\tilde{j}}^{-2d_k}}{b'\mathbf{A}\mathbf{f}(\omega_{\tilde{j}})\mathbf{A}'b}=O_p\left(\frac{\omega_{\tilde{j}}^{-2d_k}(j^{d_k-d_{k-1}}+\omega_{\tilde{j}}^{d_k-d_{k+1}}+\omega_{\tilde{j}}^{\rho})}{\omega_{\tilde{j}}^{-2d_k}}\right)=O_p(j^{d_k-d_{k-1}}).$$

Thus, the first term of (25) is

$$\frac{1}{m_n}\sum_{j=1}^{m_n}(\alpha_j-1)\left(1-\frac{\mathcal{G}\omega_{\tilde{j}}^{-2d_k}}{b'\mathbf{Af}(\omega_{\tilde{j}})\mathbf{A}'b}\right)\frac{I_{vv}(\omega_j)}{\mathcal{G}\omega_{\tilde{j}}^{-2d_k}} = O_p\left(\frac{1}{m_n}\sum_{j=1}^{m_n}(\alpha_j+1)j^{d_k-d_{k-1}}\right)$$

$$(28) \qquad = O_p\left(\frac{1}{m_n}\left(\sum_{j=1}^{m_n}\alpha_j^2\right)^{1/2}+m_n^{d_k-d_{k-1}}\right)$$

$$= o_p(1)$$

since $\sum_{j=1}^{m_n}\alpha_j^2 = O(m_n^{4(d_k-\Delta)}+m\log m)$, by Equation 3.24 of Robinson [17]. Applying (27) and (26), the second term of (25) is

$$O_p\left(\frac{1}{m_n}\sum_{j=1}^{m_n}(\alpha_j+1)j^{-\rho/2}\right) = O_p\left(\frac{1}{m_n}\left(\sum_{j=1}^{m_n}\alpha_j^2\right)^{1/2}+m_n^{-\rho/2}\right) = o_p(1),$$

by the same argument as for (28). The third term of (25) is bounded by

$$\left|\frac{1}{m_n}\sum_{j=1}^{[e^v]}(\alpha_j-1)\frac{b'\mathbf{AS}(\omega_j)\mathbf{A}'b}{b'\mathbf{Af}(\omega_{\tilde{j}})\mathbf{A}'b}\right| + \left|\frac{1}{m_n}\sum_{j=[e^v]+1}^{m_n}(\alpha_j-1)\frac{b'\mathbf{AS}(\omega_j)\mathbf{A}'b}{b'\mathbf{Af}(\omega_{\tilde{j}})\mathbf{A}'b}\right|.$$

Following from (24) and the lower bound of $b'\mathbf{Af}(\omega_{\tilde{j}})\mathbf{A}'b$ in Lemma 20, the first term of the above equation is

$$O_p\left(\omega_{m_n}^{2(d_k-\Delta)}\frac{1}{m_n}\sum_{j=1}^{[e^v]}\omega_{\tilde{j}}^{2\Delta}b'\mathbf{AS}(\omega_j)\mathbf{A}'b\right) = O_p(\omega_{m_n}^{2(d_k-\Delta)}\mathcal{M}_{m_n}(\Delta))$$

$$= o_p(\omega_{m_n}^{2(d_k-\Delta)}\omega_{m_n}^{2(\Delta-d_k)}) = o_p(1),$$

by (ii) of Lemma 21, because $0\geq \Delta-d_k > -1/2$. We will complete the proof by showing that

$$(29) \qquad \frac{1}{m_n}\sum_{j=[e^v]+1}^{m_n}\left(\frac{j}{m_n}\right)^{2(\Delta_1-d_k)}\frac{b'\mathbf{AS}(\omega_j)\mathbf{A}'b}{b'\mathbf{Af}(\omega_{\tilde{j}})\mathbf{A}'b} = o_p(1).$$

Note that $e^\nu \sim m_n/e$. Following the similar computation for (55),

$$\mathbb{E}\left\|\frac{1}{m_n}\sum_{j=[e^v]+1}^{m_n}\left(\frac{j}{m_n}\right)^{2(\Delta_1-d_k)}\frac{\mathbf{S}_{h\ell}(\omega_j)}{\omega_j^{-2d_k}}\right\|^2$$

$$= O\left(n^{-4d_k+2d_h+2d_\ell}m_n^{2(2d_k-2\Delta_1-1)}\sum_{j=\delta m_n}^{m_n}j^{4\Delta_1-2d_h-2d_\ell}\right)$$

$$= O(\omega_{m_n}^{4d_k-2d_h-2d_\ell}m_n^{-1}).$$

Hence,

$$\left\|\frac{1}{m_n}\sum_{j=[e^v]+1}^{m_n}\left(\frac{j}{m_n}\right)^{2(\Delta_1-d_k)}\frac{\mathbf{S}_{h\ell}(\omega_j)}{\omega_j^{-2d_k}}\right\| = O_p(\omega_{m_n}^{2d_k-d_h-d_\ell}m_n^{-1/2}).$$

By Lemma 16, we have

$$\begin{aligned}
&b'\mathbf{A}_h\left(\frac{1}{m_n}\sum_{j=[e^v]+1}^{m_n}\left(\frac{j}{m_n}\right)^{2(\Delta_1-d_k)}\frac{\mathbf{S}_{h\ell}(\omega_j)}{\omega_j^{-2d_k}}\right)\mathbf{A}'_\ell b\\
&\quad=\begin{cases} O_p(m_n^{2d_k-d_h-d_\ell-1/2}), & h,\ell<k,\\ O_p(\omega_{m_n}^{2d_k-d_h-d_\ell}m_n^{-1/2}), & h,\ell\geq k,\\ O_p(\omega_{m_n}^{d_k-d_\ell}m_n^{-1/2+d_k-d_h}), & h<k,\ell\geq k,\end{cases}\\
&\quad=o_p(1).
\end{aligned}$$

Equation (29) follows from the triangle inequality. □

8.3.2. *Proof of Theorem 3.* By Theorem 2, $\widehat{d}_k$ satisfies

$$0=\frac{\partial R(\hat{d}_k)}{\partial d}=\frac{\partial R(d_k)}{\partial d}+\frac{\partial^2 R(\tilde{d})}{\partial d^2}(\hat{d}_k-d_k), \tag{30}$$

where $|\tilde{d}-d_k|\leq|\hat{d}_k-d_k|$. Let

$$\mathbf{Z}_n=2m_n^{-1/2}\sum_{j=1}^{m_n}\nu_j(I_{\varepsilon\varepsilon}(\omega_j)-\mathbf{\Sigma}),\qquad \nu_j=\log\tilde{j}-\frac{1}{m_n}\sum_{j=1}^{m_n}\log\tilde{j}$$

and

$$\mathfrak{Z}_n=\frac{1}{\mathcal{G}}b'\mathbf{A}_k\mathbf{\Psi}_k^{\dagger\prime}(0)\mathbf{Z}_n\mathbf{\Psi}_k^{\dagger}(0)\mathbf{A}'_k b,$$

where

$$\mathcal{G}=b'\mathbf{A}_k\mathbf{f}^{\dagger}(0)\mathbf{A}'_k b=b'\mathbf{A}_k\mathbf{\Psi}_k^{\dagger\prime}(0)\mathbf{\Sigma}\mathbf{\Psi}_k^{\dagger}(0)\mathbf{A}'_k b \tag{31}$$

and $\mathbf{\Psi}_k^{\dagger}(\omega)$ is a $q\times a_k$ submatrix of $\mathbf{\Psi}^{\dagger}(\omega)=[\mathbf{\Psi}_0^{\dagger}(\omega)\ \ldots\ \mathbf{\Psi}_s^{\dagger}(\omega)]$ in (8). We show in Lemmas 13 and 14 that

$$\frac{\partial^2 R(\tilde{d})}{\partial d^2}\xrightarrow{p}4 \tag{32}$$

and

$$m_n^{1/2}\frac{\partial R(d_k)}{\partial d}=\mathfrak{Z}_n+o_p(1). \tag{33}$$

From Lemmas 0 and 8 of Hurvich and Chen [8], the $(u,v)$th entry of $\mathbf{Z}_n$ satisfies

$$Z_{n,uv} \xrightarrow{D} N(0, 4\Phi_p \sigma_{uv}^2).$$

Using a similar computation for the variance above and Eq. (53) in the proof of Lemma 19, we obtain

$$E(Z_{n,u_1v_1} Z_{n,u_2v_2}) \to 4\Phi_p \sigma_{u_1v_2}\sigma_{u_2v_1}.$$

Using the Cramer–Wold device, we have $\operatorname{vec}\mathbf{Z}_n \xrightarrow{D} \operatorname{vec}\mathbf{Z} \sim N(0, 4\Phi_p \mathbf{\Sigma} \otimes \mathbf{\Sigma})$. By Lemma 15, $b \xrightarrow{D} \mathring{b}$, thus,

$$\mathfrak{z}_n \xrightarrow{D} \frac{\mathring{b}'\mathbf{A}_k \mathbf{\Psi}_k^{\dagger\prime}(0)\mathbf{Z}\mathbf{\Psi}_k^{\dagger}(0)\mathbf{A}_k'\mathring{b}}{\mathring{b}'\mathbf{A}_k \mathbf{\Psi}_k^{\dagger\prime}(0)\mathbf{\Sigma}\mathbf{\Psi}_k^{\dagger}(0)\mathbf{A}_k'\mathring{b}} := \mathfrak{z}. \tag{34}$$

Let $\varphi = (\varphi_1, \ldots, \varphi_q)' = \mathbf{\Psi}_k^{\dagger\prime}(0)\mathbf{A}_k'\mathring{b}$. Then $\mathfrak{z}|\mathring{b}$ is a normal random variable with mean zero and variance

$$\operatorname{var}(\mathfrak{z}|\mathring{b}) = 4\Phi_p \frac{\varphi'\mathbf{\Sigma}\varphi\varphi'\mathbf{\Sigma}\varphi}{(\varphi'\mathbf{\Sigma}\varphi)^2} = 4\Phi_p.$$

Thus, $\mathfrak{z}$ is independent of $\mathring{b}$ and $\mathfrak{z} \sim N(0, 4\Phi_p)$. Together with (30) and (32)–(34), we have proved the theorem. □

LEMMA 13. *Let $\tilde{d}$ be such that $|\tilde{d} - d_k| \le |\hat{d}_k - d_k|$. Then under the assumptions of Theorem 2,*

$$\frac{\partial^2 R(\tilde{d})}{\partial d^2} \xrightarrow{p} 4.$$

PROOF. Define

$$\widehat{G}_a(d) = \frac{1}{m_n}\sum_{j=1}^{m_n} (\log \omega_{\tilde{j}})^a \omega_{\tilde{j}}^{2d} I_{vv}(\omega_j)$$

and

$$\widehat{F}_a(d) = \frac{1}{m_n}\sum_{j=1}^{m_n} (\log \tilde{j})^a \omega_{\tilde{j}}^{2d} I_{vv}(\omega_j), \qquad \widehat{E}_a(d) = \frac{1}{m_n}\sum_{j=1}^{m_n} (\log \tilde{j})^a \tilde{j}^{2d} I_{vv}(\omega_j).$$

Then

$$\begin{aligned} \frac{\partial^2 R(d)}{\partial d^2} &= \frac{4\{\widehat{G}_2(d)\widehat{G}(d) - \widehat{G}_1^2(d)\}}{\widehat{G}^2(d)} = \frac{4\{\widehat{F}_2(d)\widehat{F}_0(d) - \widehat{F}_1^2(d)\}}{\widehat{F}_0^2(d)} \\ &= \frac{4\{\widehat{E}_2(d)\widehat{E}_0(d) - \widehat{E}_1^2(d)\}}{\widehat{E}_0^2(d)}. \end{aligned} \tag{35}$$

We first show that

$$\widehat{F}_a(\tilde{d}) = \widehat{F}_a(d_k) + o_p(1), \qquad a = 1, 2, 3, \tag{36}$$

by showing that $\widehat{E}_a(\tilde{d}) = \widehat{E}_a(d_k) + o_p(n^{2d_k})$ for $a = 0, 1, 2$. Let $M = \{d : \log^3 m_n \times |d - d_k| \le \varepsilon\}$, where $\varepsilon > 0$ is fixed to be such that $2\varepsilon < \log^2 m_n$ with a proper $n$. Following the same line of proof as on page 1642 of Robinson [17], for $\eta > 0$,

$$\begin{aligned} &P\biggl(|\widehat{E}_a(\tilde{d}) - \widehat{E}_a(d_k)| > \eta\biggl(\frac{2\pi}{n}\biggr)^{-2d_k}\biggr) \\ &\quad \le P\biggl(\widehat{G}(d_k) > \tfrac{\eta}{2e\varepsilon}(\log m_n)^{2-a}\biggr) + P(\log^3 m_n |\tilde{d} - d_k| > \varepsilon). \end{aligned} \tag{37}$$

The first probability is bounded by

$$P\biggl(|\widehat{G}(d_k) - \mathcal{G}| > \frac{\eta}{4e\varepsilon}(\log m_n)^{2-a}\biggr) + P\biggl(\mathcal{G} > \frac{\eta}{4e\varepsilon}(\log m_n)^{2-a}\biggr).$$

Both probabilities in the above equation tend to 0 for $\varepsilon$ sufficiently small since $|\widehat{G}(d_k) - \mathcal{G}| = o_p(1)$ and $\mathcal{G} < C$, by Lemma 20. To show that the second probability in (37) tends to 0, we only have to verify that

$$\sup_{\Theta_1 \cap \mathcal{N}_\delta} \left| \frac{\widehat{G}(d) - G(d)}{G(d)} \right| = o_p(\log^{-6} m_n).$$

From (22) and (23) in the proof of Lemma 11,

$$\sup_{\Theta_1 \cap \mathcal{N}_\delta} \left| \frac{\widehat{G}(d) - G(d)}{G(d)} \right| \le \sup_{\Theta_1} \left| \frac{\widehat{G}(d) - G(d)}{G(d)} \right| = o_p(m_n^{-\varepsilon}).$$

We have established (36). Combining this with (35), we have

$$\frac{\partial^2 R(\tilde{d})}{\partial d^2} = \frac{4\{\widehat{F}_2(d_k)\widehat{F}_0(d_k) - \widehat{F}_1^2(d_k)\}}{\widehat{F}_0^2(d_k)} + o_p(1) \qquad \text{as } n \to \infty.$$

By Lemma 21,

$$\begin{aligned} \left| \widehat{F}_a(d_k) - \mathcal{G}\frac{1}{m_n}\sum_{j=1}^{m_n} \log^a \tilde{j} \right| &= \left| \frac{1}{m_n}\sum_{j=1}^{m_n} \log^a \tilde{j} \left( \frac{I_{vv}(\omega_j)}{\omega_{\tilde{j}}^{-2d_k}} - \mathcal{G} \right) \right| \\ &\le \log^a m_n |\mathcal{L}_{m_n}(d_k) + \mathcal{M}_{m_n}(d_k) + \mathcal{F}_{m_n}(d_k) - \mathcal{G}| \\ &= O_p(m_n^{-\varepsilon} \log^a m_n). \end{aligned}$$

By the same reasoning as that used in (4.10) of Robinson [17], we obtain

$$\frac{\partial^2 R(\tilde{d})}{\partial d^2} = 4\left\{ \frac{1}{m_n}\sum_{j=1}^{m_n} \log^2 \tilde{j} - \left( \frac{1}{m_n}\sum_{j=1}^{m_n} \log \tilde{j} \right)^2 \right\}(1 + o_p(1)) + o_p(1) \xrightarrow{p} 4. \qquad \square$$

LEMMA 14. *Under the assumptions of Theorem 3,*

$$m_n^{1/2}\frac{\partial R(d_k)}{\partial d}-\frac{1}{\mathcal{G}}b'\mathbf{A}_k\mathbf{\Psi}_k^{\dagger *}(0)\mathbf{Z}_n\mathbf{\Psi}_k^{\dagger}(0)\mathbf{A}_k'b=o_p(1),$$

*where*

$$\mathbf{Z}_n=2m_n^{-1/2}\sum_{j=1}^{m_n}\nu_j(I_{\varepsilon\varepsilon}(\omega_j)-\mathbf{\Sigma})\quad \text{and}\quad \nu_j=\log\tilde{j}-\frac{1}{m_n}\sum_{j=1}^{m_n}\log\tilde{j}.$$

PROOF. Note that

$$\frac{\partial R(d)}{\partial d}=\frac{2}{m_n}\sum_{j=1}^{m_n}\frac{\upsilon_j I_{vv}(\omega_j)}{\omega_{\tilde{j}}^{-2d}\widehat{G}(d)}.$$

Since $\widehat{G}(d_k)-\mathcal{G}=o_p(1)$, by (23), and $\sum_{j=1}^{m_n}\nu_j=0$, we have

$$m_n^{1/2}\frac{\partial R(d_k)}{\partial d}=2m_n^{-1/2}\sum_{j=1}^{m_n}\nu_j\left(\frac{I_{vv}(\omega_j)}{\mathcal{G}\omega_{\tilde{j}}^{-2d_k}}-1\right)(1+o_p(1))$$

and

$$\begin{aligned}
&m_n^{1/2}\frac{\partial R(d_k)}{\partial d}-\frac{1}{\mathcal{G}}b'\mathbf{A}_k\mathbf{\Psi}_k^{\dagger *}(0)\mathbf{Z}_n\mathbf{\Psi}_k^{\dagger}(0)\mathbf{A}_k'b\\
&\quad=2m_n^{-1/2}\left\{\sum_{j=1}^{m_n}\nu_j\left(\frac{I_{vv}(\omega_j)}{\mathcal{G}\omega_{\tilde{j}}^{-2d_k}}-1\right)-\sum_{j=1}^{m_n}\frac{\nu_j}{\mathcal{G}\omega_{\tilde{j}}^{-2d_k}}b'\mathbf{A}_k\mathbf{S}_{kk}(\omega_j)\mathbf{A}_k'b\right\}\\
&\qquad+2m_n^{-1/2}\left\{\sum_{j=1}^{m_n}\frac{\nu_j}{\mathcal{G}}(\omega_{\tilde{j}}^{2d_k}b'\mathbf{A}_k\mathbf{S}_{kk}(\omega_j)\mathbf{A}_k'b)-b'\mathbf{A}_k\mathbf{\Psi}_k^{\dagger *}(0)\mathbf{Z}_n\mathbf{\Psi}_k^{\dagger}(0)\mathbf{A}_k'b\right\} \qquad (38)\\
&\qquad+o_p(1),
\end{aligned}$$

where $\mathbf{S}_{h\ell}(\omega_j)$ is the $(h,\ell)$th block of $\mathbf{S}(\omega_j)$ defined in (20). Let

$$\mathcal{M}_{m_n}^{(h,\ell)}(d)=\frac{1}{m_n}\sum_{j=1}^{m_n}\omega_{\tilde{j}}^{2d}b'\mathbf{A}_h\mathbf{S}_{h\ell}(\omega_j)\mathbf{A}_\ell'b.$$

The first term of (38) is then

$$\begin{aligned}
&\left|2m_n^{-1/2}\sum_{j=1}^{m_n}\frac{\nu_j}{\mathcal{G}\omega_{\tilde{j}}^{-2d_k}}\left\{I_{vv}(\omega_j)-\mathcal{G}\omega_{\tilde{j}}^{-2d_k}+\sum_{\substack{u,v=0\\ v\neq k}}^{q}b'\mathbf{A}_h\mathbf{S}_{h\ell}(\omega_j)\mathbf{A}_\ell'b\right\}\right|\\
&\quad\le\frac{2m_n^{1/2}\log m_n}{\mathcal{G}}\left|\mathcal{L}_{m_n}(d_k)+\mathcal{F}_{m_n}(d_k)-\mathcal{G}+\sum_{\substack{u,v=0\\ v\neq k}}^{q}\mathcal{M}_{m_n}^{(h,\ell)}(d_k)\right|\\
&\quad=o_p(1),
\end{aligned}$$

by (20), (21), Lemma 22 and Assumption 3B. Since

$$\begin{aligned}\mathbf{S}_{kk}(\omega_j) &= |1-e^{i\omega_j}|^{-2d_k}\mathbf{\Psi}_k^{\dagger *}(\omega_j)[I(\omega_j)-\mathbf{\Sigma}]\mathbf{\Psi}_k^{\dagger}(\omega_j)\\ &= \omega_{\tilde{j}}^{-2d_k}\mathbf{\Psi}_k^{\dagger *}(0)[I(\omega_j)-\mathbf{\Sigma}]\mathbf{\Psi}_k^{\dagger}(0)+O_p(\omega_{\tilde{j}}^{-2d_k+\rho}),\end{aligned}$$

the second term in (38) is

$$O_p\left(m_n^{-1/2}\sum_{j=1}^{m_n}\nu_j\omega_{\tilde{j}}^{\rho}\right)=O_p\left(\frac{m_n^{\rho+1/2}\log m}{n^{\rho}}\right)=o_p(1),$$

by Lemma 16 and Assumption 3B. We have shown that both terms on the right-hand side of (38) are $o_p(1)$ and, hence, have completed the proof. □

LEMMA 15. *Under Assumption 1, the matrix* $\mathbf{X}=\mathbf{X}\{I_m(y_t,y_t)\}$ *satisfies*

$$\mathbf{X}\xrightarrow{D}\mathring{\mathbf{X}}=\mathbf{B}'\mathbf{X}\{\mathbf{K}(\mathbf{U}\mathbf{U}'+\mathbf{V}\mathbf{V}')_D\mathbf{K}'\},$$

*where* $\mathbf{X}\{\mathbf{K}(\mathbf{U}\mathbf{U}'+\mathbf{V}\mathbf{V}')_D\mathbf{K}'\}$ *is the matrix of normalized eigenvectors of* $\mathbf{K}(\mathbf{U}\mathbf{U}'+\mathbf{V}\mathbf{V}')_D\mathbf{K}'$ *in Lemma 6 and* $\mathbf{U},\mathbf{V}$ *are defined as in Lemma 1. Thus,* $\mathring{\mathbf{X}}$ *is a continuous function with respect to* $\operatorname{vec}(\mathbf{U},\mathbf{V})$.

PROOF. It suffices to show that

$$\|\mathbf{X}_k(\mathbf{H})-\mathbf{X}_k\|_F=O_p(n^{-\alpha_k}) \tag{39}$$

and that the eigenvectors of $\mathbf{H}$ satisfy

$$\chi_j(\mathbf{H})\xrightarrow{D}\mathring{\xi}_j(\operatorname{vec}(\mathbf{U},\mathbf{V})),\qquad j=1,\dots,q, \tag{40}$$

where $\mathring{\xi}_j$ are continuous functions of $\operatorname{vec}(\mathbf{U},\mathbf{V})$.

We first show (39). Note that since both $\mathbf{H}$ and $I_m$ are symmetric, we can assume that $\chi_j'\chi_j(\mathbf{H})\ge 0$. We have

$$\begin{aligned}\|\mathbf{X}_k(\mathbf{H})-\mathbf{X}_k\|_F^2 &\le a_k\max_{j\in N_k}\|\chi_j(\mathbf{H})-\chi_j\|^2\le C\max_{j\in N_k}\sin^2\theta(\chi_j,\chi_j(\mathbf{H}))\\ &\le C\|\sin\Theta\{\mathcal{M}(\mathbf{X}_k(\mathbf{H})),\mathcal{M}(\mathbf{X}_k)\}\|_F^2,\end{aligned}$$

by the definition of the $\sin\Theta$ bound. Equation (39) follows from Lemma 4.

Next, we derive (40). Since $\chi_j(\mathbf{H})=\mathbf{B}'\chi_j(\mathbf{\Phi}_D)$, it is sufficient to show that

$$\chi_j(\mathbf{\Phi}_D)\xrightarrow{D}\mathring{\varsigma}_j(\operatorname{vec}(\mathbf{U},\mathbf{V})),\qquad j=1,\dots,q,$$

where $\mathring{\varsigma}_j$ are continuous functions of $\operatorname{vec}(\mathbf{U},\mathbf{V})$. Let $\widetilde{\mathbf{\Phi}}_D=\mathbf{d}_n^{-1}\mathbf{\Phi}_D\mathbf{d}_n^{-1}$. Then

$$\widetilde{\mathbf{\Phi}}_D=\mathbf{X}'(\widetilde{\mathbf{\Phi}}_D)\mathbf{\Lambda}(\widetilde{\mathbf{\Phi}}_D)\mathbf{X}(\widetilde{\mathbf{\Phi}}_D).$$

First, note that the eigenvalues of $\widetilde{\mathbf{\Phi}}_D$ are distinct, with probability 1, by Okamoto [14]. Since both $\widetilde{\mathbf{\Phi}}_D$ and $\mathbf{\Phi}_D$ are block diagonal matrices, we have

$$\mathbf{\Phi}_D = \mathbf{d}_n^{-1}\mathbf{X}'(\widetilde{\mathbf{\Phi}}_D)\mathbf{\Lambda}(\widetilde{\mathbf{\Phi}}_D)\mathbf{X}(\widetilde{\mathbf{\Phi}}_D)\mathbf{d}_n^{-1} = \mathbf{X}'(\widetilde{\mathbf{\Phi}}_D)\mathbf{d}_n^{-1}\mathbf{\Lambda}(\widetilde{\mathbf{\Phi}}_D)\mathbf{d}_n^{-1}\mathbf{X}(\widetilde{\mathbf{\Phi}}_D).$$

This implies that

$$\mathbf{X}'(\mathbf{\Phi}_D) = \mathbf{X}'(\widetilde{\mathbf{\Phi}}_D) \quad \text{and} \quad \mathbf{\Lambda}(\mathbf{\Phi}_D) = \mathbf{d}_n^{-1}\mathbf{\Lambda}(\widetilde{\mathbf{\Phi}}_D)\mathbf{d}_n^{-1}.$$

We now let $\mathbf{K}$ be defined as in Lemma 6 and rewrite $\mathbf{U}_n$ and $\mathbf{V}_n$ in (10) as

$$\mathbf{U}_n = [\mathbf{U}'_{n,0} \ \ldots \ \mathbf{U}'_{n,s}]', \qquad \mathbf{V}_n = [\mathbf{V}'_{n,0} \ \ldots \ \mathbf{V}'_{n,s}]',$$

where $\mathbf{U}_{n,k}$ and $\mathbf{V}_{n,k}$ are $a_k \times m$ matrices. Since $\mathbf{K}$ is a block diagonal matrix, we have

$$\widetilde{\mathbf{\Phi}}_D = \mathbf{K}' \operatorname{diag}(\mathbf{U}_{n,0}\mathbf{U}'_{n,0} + \mathbf{U}_{n,0}\mathbf{U}'_{n,0}, \ldots, \mathbf{U}_{n,s}\mathbf{U}'_{n,s} + \mathbf{U}_{n,s}\mathbf{U}'_{n,s})\mathbf{K}.$$

It follows that

$$\chi_j(\mathbf{\Phi}_D) = \chi_j(\widetilde{\mathbf{\Phi}}_D) := \mathring{\varsigma}_j(\operatorname{vec}(\mathbf{U}_n, \mathbf{V}_n)) \xrightarrow{D} \mathring{\varsigma}_j(\operatorname{vec}(\mathbf{U}, \mathbf{V})),$$

where $\mathring{\varsigma}_j(\cdot)$ is a continuous function because the eigenvalues of $\widetilde{\mathbf{\Phi}}_D$ are distinct with probability 1 and

$$\chi_j(\mathbf{H}) = \mathbf{B}'\chi_j(\mathbf{\Phi}_D) \xrightarrow{D} \mathbf{B}'\mathring{\varsigma}_j(\operatorname{vec}(\mathbf{U}, \mathbf{V})) = \mathring{\xi}_j(\operatorname{vec}(\mathbf{U}, \mathbf{V})). \qquad \square$$

REMARK 3. Let $\mathring{\mathbf{X}} = [\mathring{\mathbf{X}}_0 \ \ldots \ \mathring{\mathbf{X}}_s]$. Since $\mathbf{X}\{\mathbf{K}(\mathbf{U}\mathbf{U}' + \mathbf{V}\mathbf{V}')_D\mathbf{K}'\}$ is a block diagonal matrix, $\mathring{\mathbf{X}}_k = \mathbf{B}'\mathbf{X}_k\{\mathbf{K}(\mathbf{U}\mathbf{U}' + \mathbf{V}\mathbf{V}')_D\mathbf{K}'\} \in \mathcal{B}_k$.

8.3.3. *Proof of Theorem 4.* In case of no cointegration, we have

$$C_1 = \lambda_1(\mathbf{A}\mathbf{A}') \geq \|b'\mathbf{A}\|^2 = (b'\mathbf{A}\mathbf{A}'b) \geq \lambda_q(\mathbf{A}\mathbf{A}') \geq C_2$$

and

$$\tilde{C}_1 = \lambda_1(\mathbf{A}\mathbf{A}')\lambda_1(\mathbf{f}^\dagger(0)) \geq \mathcal{G} = b'\mathbf{A}\mathbf{f}^\dagger(0)\mathbf{A}'b \geq \lambda_q(\mathbf{A}\mathbf{A}')\lambda_q(\mathbf{f}^\dagger(0)) = \tilde{C}_1,$$

where $C_1$, $C_2$, $\tilde{C}_1$ and $\tilde{C}_2$ are positive constants. Furthermore, by Assumption 2,

$$b'\mathbf{A}\mathbf{f}(\omega_j)\mathbf{A}'b - \omega_j^{-2d_0}\mathcal{G} = O_p(\omega_j^{-2d_0+\rho}),$$

for $1 \leq j \leq m_n$. Following along the lines of the proofs of Theorems 2 and 3, we have $m_n^{1/2}(\hat{d}_{aa} - d_{aa}) \xrightarrow{D} N(0, \Phi_p/4)$ for $a = 1, \ldots, q$. The theorem follows. $\square$

**9. Technical lemmas.** We will need the following two lemmas:

LEMMA 16. *If $b = \chi_a$, where $a \in N_k$, then, under Assumption 1,*

$$b'\mathbf{A}_h = O_p(n^{-d_h+d_k})$$

*for $h < k$, $k = 1, \dots, s$ and*

$$b'\mathbf{A}_k = O_p(1),$$

*for $k = 0, \dots, s$.*

LEMMA 17. *If $b = \chi_a$, where $a \in N_k$, then, under Assumption 1,*

$$\|b'\mathbf{A}_k\| \ge C(1 - \varepsilon_k),$$

*where $C > 0$ and $\varepsilon_k = O_p(n^{-\alpha_k})$, $k = 0, \dots, s$.*

PROOF OF LEMMA 16. Since $\mathbf{X}(\mathbf{H})$ is an orthogonal matrix and $\mathcal{M}\mathbf{X}(\mathbf{H}) = \mathbb{R}^q$, we have

$$b = \sum_{\ell=0}^{s} \mathbf{X}_\ell(\mathbf{H}) c_\ell, \tag{41}$$

where

$$c_\ell = \mathbf{X}'_\ell(\mathbf{H}) b = O_p(n^{-|d_k - d_\ell|}),$$

by Lemma 7. Furthermore, for $\ell > h$,

$$\begin{aligned}
\mathbb{E}[\|\mathbf{X}'_\ell(\mathbf{H})\mathbf{A}_h\|] &= \mathbb{E}[\|\mathbf{X}'_\ell(\mathbf{H})\mathbf{A}_h\| \mathbf{1}_{\{\mathcal{M}\mathbf{X}_\ell(\mathbf{H}) \subset \oplus_{j>h} \mathcal{B}_j\}}] \\
&\quad + \mathbb{E}[\|\mathbf{X}'_\ell(\mathbf{H})\mathbf{A}_h\| \mathbf{1}_{\{\mathcal{M}\mathbf{X}_\ell(\mathbf{H}) \cap \oplus_{j\le h} \mathcal{B}_j \ne \mathbf{0}\}}] \\
&\le 0 + \mathbb{E}[\|\mathbf{X}'_\ell(\mathbf{H})\mathbf{A}_h\| \mathbf{1}_{\{\mathcal{M}\mathbf{X}_\ell(\mathbf{H}) \cap \oplus_{j\le h} \mathcal{B}_j \ne \mathbf{0}\}}] \\
&= \mathbb{E}[\operatorname{trace}^{1/2}(\mathbf{A}'_h \mathbf{X}_\ell(\mathbf{H}) \mathbf{X}'_\ell(\mathbf{H}) \mathbf{A}_h) \mathbf{1}_{\{\mathcal{M}\mathbf{X}_\ell(\mathbf{H}) \cap \oplus_{j\le h} \mathcal{B}_j \ne \mathbf{0}\}}] \\
&\le \mathbb{E}[\operatorname{trace}^{1/2}(\mathbf{A}'_h \mathbf{A}_h) \operatorname{trace}^{1/2}(\mathbf{X}_\ell(\mathbf{H}) \mathbf{X}'_\ell(\mathbf{H})) \mathbf{1}_{\{\mathcal{M}\mathbf{X}_\ell(\mathbf{H}) \cap \oplus_{j\le h} \mathcal{B}_j \ne \mathbf{0}\}}] \\
&= a_\ell^{1/2} \|\mathbf{A}_h\| \left| P\left\{ \mathcal{M}\mathbf{X}_\ell(\mathbf{H}) \cap \bigoplus_{j\le h} \mathcal{B}_j \ne \mathbf{0} \right\} \right|^{1/2} \\
&= O(n^{-d_h + d_\ell}),
\end{aligned} \tag{42}$$

by Lemma 10. For $\ell \le h$,

$$\mathbb{E}[\|\mathbf{X}'_\ell(\mathbf{H})\mathbf{A}_h\|] = O(1). \tag{43}$$

We have, for $h < k$,

$$\begin{aligned} b'\mathbf{A}_h &= \sum_{\ell=0}^{s} c_\ell' \mathbf{X}_\ell'(\mathbf{H})\mathbf{A}_h \\ &= \sum_{\ell:\ell\le h} c_\ell' \mathbf{X}_\ell'(\mathbf{H})\mathbf{A}_h + \sum_{\ell:\ell>h} c_\ell' \mathbf{X}_\ell'(\mathbf{H})\mathbf{A}_h \\ &= O_p\biggl(\sum_{\ell:\ell\le h} n^{-d_\ell+d_k} + \sum_{\ell:h<\ell\le k} n^{-d_h+d_\ell-d_\ell+d_k} + \sum_{\ell:\ell>k} n^{-d_h+d_\ell-d_k+d_\ell}\biggr) \\ &= O_p(n^{-d_h+d_k}). \end{aligned}$$

For $h = k$, the above equation is of $O_p(1)$ since $c_k = O_p(1)$ and $\mathbb{E}[\|\mathbf{X}_k'(\mathbf{H})\mathbf{A}_k\|] = O(1)$. $\square$

Proof of Lemma 17. Note that

$$\begin{aligned} \|b'\mathbf{A}_k\|^2 &= \left\| c_k' \mathbf{X}_k'(\mathbf{H})\mathbf{A}_k + \sum_{\ell=0,\ell\ne k}^{s} c_\ell' \mathbf{X}_\ell'(\mathbf{H})\mathbf{A}_k \right\| \\ &\ge \left| \|c_k' \mathbf{X}_k'(\mathbf{H})\mathbf{A}_k\| - \left\| \sum_{\ell=0,\ell\ne k}^{s} c_\ell' \mathbf{X}_\ell'(\mathbf{H})\mathbf{A}_k \right\| \right|. \end{aligned} \tag{44}$$

Using (41), we have

$$1 = \|b\|^2 = \sum_{\ell=0}^{s} \|\mathbf{X}_\ell(\mathbf{H})c_\ell\|^2 = \sum_{\ell=0}^{s} \|c_\ell\|^2 = \|c_k\|^2 + \sum_{\ell=0,\ell\ne k}^{s} \|c_\ell\|^2$$

and

$$\sum_{\ell=0,\ell\ne k}^{s} \|c_\ell\|^2 = O_p(n^{-2\alpha_k}), \tag{45}$$

by Lemma 4. Thus,

$$\|c_k\|^2 = 1 - O_p(n^{-2\alpha_k}). \tag{46}$$

By (45), (42) and (43),

$$\begin{aligned} \left\| \sum_{\ell=0,\ell\ne k}^{s} c_\ell' \mathbf{X}_\ell'(\mathbf{H})\mathbf{A}_k \right\| &\le \left( \sum_{\ell=0,\ell\ne k}^{s} \|c_\ell\|^2 \sum_{\ell=0,\ell\ne k}^{s} \|\mathbf{X}_\ell'(\mathbf{H})\mathbf{A}_k\|^2 \right)^{1/2} \\ &= O_p(n^{-\alpha_k}). \end{aligned} \tag{47}$$

Furthermore, if $\mathcal{M}\mathbf{X}_k(\mathbf{H}) = \mathcal{B}_k$, then there exists an $a_k \times a_k$ orthogonal matrix $\mathbf{D}$ such that

$$\mathbf{X}_k(\mathbf{H}) = \mathbf{B}_k\mathbf{D}$$

since both $\mathbf{X}_k(\mathbf{H})$ and $\mathbf{B}_k$ are matrices with orthonormal columns. We have

$$\begin{aligned}\|c_k\|^2 &= \operatorname{trace}\{c_k'\mathbf{D}\mathbf{B}_k'\mathbf{A}_k(\mathbf{B}_k'\mathbf{A}_k)^{-1}(\mathbf{A}_k'\mathbf{B}_k)^{-1}\mathbf{A}_k'\mathbf{B}_k\mathbf{D}'c_k\}\\ &\le \|(\mathbf{B}_k'\mathbf{A}_k)^{-1}\|^2\|c_k'\mathbf{D}\mathbf{B}_k'\mathbf{A}_k\|^2\\ &= \|(\mathbf{B}_k'\mathbf{A}_k)^{-1}\|^2\|c_k'\mathbf{X}_k'(\mathbf{H})\mathbf{A}_k\|^2.\end{aligned}$$

It follows that

$$\|c_k'\mathbf{X}_k'(\mathbf{H})\mathbf{A}_k\|^2 \ge \|(\mathbf{B}_k'\mathbf{A}_k)^{-1}\|^{-2}\|c_k\|^2 = C(1-\delta_k),$$

where $\delta_k = O_p(n^{-2\alpha_k})$, by (46). By (44), (47) and the above equation, $\|b'\mathbf{A}_k\| \ge C(1-\delta_k-\tilde{\varepsilon}_k)$, where $\tilde{\varepsilon}_k = O_p(n^{-\alpha_k})$. We have thus completed the proof. □

LEMMA 18. *Let $R_{ab}(\omega_j)$ be the $(a,b)$th entry of $\mathbf{R}(\omega_j)$,*

$$\mathbb{E}|R_{ab}(\omega_j)| \le C|1-e^{-i\omega_{\tilde{j}}}|^{-(d_{aa}+d_{bb})}j^{-\rho/2},$$

$$a,b=1,\dots,q \text{ and } 1\le j\le [n/2]$$

*under Assumption 2.*

PROOF. Let $J_{z_a}(\omega_j)$ be the $j$th element of $J_z(\omega_j)$, the discrete Fourier transform of $z_t$. By (4),

$$J_{z_a}(\omega_j) = \sum_{b=1}^{q} J_{z_{a_b}}(\omega_j), \tag{48}$$

where

$$J_{z_{a_b}}(\omega_j) = \frac{1}{\sqrt{2\pi\sum|h_t^{p-1}|^2}}\sum_{t=1}^{n} h_t^{p-1}\left(\sum_{k=\infty}^{\infty}\psi_{k,ab}\varepsilon_{t-k,b}\right)e^{i\omega_j t}.$$

Hence,

$$\begin{aligned}R_{ab}(\omega_j) &= J_{z_a}(\omega_j)\overline{J}_{z_b}(\omega_j) - \sum_{u=1}^{q}\Psi_{au}(\omega_{\tilde{j}})J_{\varepsilon_u}(\omega_j)\sum_{v=1}^{q}\overline{\Psi}_{bv}(\omega_{\tilde{j}})\overline{J}_{\varepsilon_v}(\omega_j)\\ &= \sum_{u,v=1}^{q}(\Psi_{au}(\omega_{\tilde{j}})\overline{\Psi}_{bv}(\omega_{\tilde{j}})(A_{au,j}\overline{A}_{bv,j} - B_{u,j}\overline{B}_{v,j})),\end{aligned} \tag{49}$$

where

$$A_{au,j} = \frac{J_{z_{a_u}}(\omega_j)}{\Psi_{au}(\omega_{\tilde{j}})} \quad \text{and} \quad B_{u,j} = J_{\varepsilon_u}(\omega_j). \tag{50}$$

From Lemmas 9 and 10 of Hurvich et al. [9],

$$\begin{aligned}\mathbb{E}|A_{au,j}-B_{u,j}|^{2\ell} &\le C\left(\int_{-\pi}^{\pi}\left|\frac{\Psi_{au}(\omega)}{\Psi_{au}(\omega_{\tilde{j}})}-1\right|^2|D_{p,n}(\omega_j-\omega)|^2\,d\omega\right)^{\ell}\\ &\le Cj^{-\ell\rho}.\end{aligned} \tag{51}$$

By the Cauchy–Schwarz inequality,

$$
\begin{aligned}
\mathbb{E}|A_{au,j}\overline{A}_{bv,j} - B_{u,j}\overline{B}_{v,j}|^2 &= \mathbb{E}|(A_{au,j} - B_{u,j})(\overline{A}_{bv,j} - \overline{B}_{v,j}) \\
&\qquad + B_{u,j}(\overline{A}_{bv,j} - \overline{B}_{v,j}) + \overline{B}_{v,j}(A_{au,j} - B_{u,j})|^2 \\
&\le 3(\mathbb{E}|A_{au,j} - B_{u,j}|^4 \mathbb{E}|A_{bv,j} - B_{v,j}|^4)^{1/2} \\
&\qquad + (\mathbb{E}|B_{u,j}|^4 \mathbb{E}|A_{bv,j} - B_{v,j}|^4)^{1/2} \\
&\qquad + (\mathbb{E}|B_{v,j}|^4 \mathbb{E}|A_{au,j} - B_{u,j}|^4)^{1/2} \\
&\le C[(j^{-2\rho} j^{-2\rho})^{1/2} + (j^{-2\rho})^{1/2}] = Cj^{-\rho}.
\end{aligned} \tag{52}
$$

We have, from (49) and Assumption 2,

$$
\begin{aligned}
\mathbb{E}|R_{ab}| &\le \sum_{u,v=1}^{q} \Psi_{au}(\omega_{\tilde{j}})\overline{\Psi}_{bv}(\omega_{\tilde{j}})(\mathbb{E}|A_{au,j}\overline{A}_{bv,j} - B_{u,j}\overline{B}_{v,j}|^2)^{1/2} \\
&\le C \sum_{u,v=1}^{q} |1 - e^{-i\omega_{\tilde{j}}}|^{-(d_{au}+d_{bv})} \tau_{au}(\omega_j)\tau_{av}(\omega_j) j^{-\rho/2} \\
&\le C|1 - e^{-i\omega_{\tilde{j}}}|^{-(d_{aa}+d_{bb})} j^{-\rho/2},
\end{aligned}
$$

where the constant $C$ does not depend on $n$. □

LEMMA 19. *Let $S_{ab}(\omega)$ be the $(a,b)$th entry of $\mathbf{S}(\omega_j)$. Then for $1 \le j, k \le [n/2]$,*

$$
\mathbb{E}|S_{ab}(\omega_j)S_{ab}(\omega_k)| \le \begin{cases} C|(1 - e^{-i\omega_{\tilde{j}}})(1 - e^{-i\omega_{\tilde{k}}})|^{-(d_{aa}+d_{bb})}, & |j-k| < p, \\ C/n, & \text{otherwise}, \end{cases}
$$

*under Assumptions 1 and 2.*

PROOF. Note that $\mathbb{E}I_{\varepsilon\varepsilon}(\omega_j) = \mathbf{\Sigma}$ and $S_{ab}(\omega_j) = \sum_{u,v=1}^{q} \Psi_{au}(\omega_{\tilde{j}})\overline{\Psi}_{bv}(\omega_{\tilde{j}}) \times (I_{\varepsilon\varepsilon,uv}(\omega_j) - \sigma_{uv})$. Now,

$$
\begin{aligned}
\mathbb{E}|S_{ab}(\omega_j)S_{ab}(\omega_k)| = \sum_{u_1,u_2,v_1,v_2=1}^{q} & \Psi_{au_1}(\omega_{\tilde{k}})\overline{\Psi}_{au_2}(\widetilde{\omega}_k)\overline{\Psi}_{bv_1}(\omega_{\tilde{j}})\Psi_{bv_2}(\omega_{\tilde{k}}) \\
& \times \mathbb{E}[(I_{\varepsilon\varepsilon,u_1v_1}(\omega_j) - \sigma_{u_1v_1})(I_{\varepsilon\varepsilon,u_2v_2}(\omega_k) - \sigma_{u_2v_2})].
\end{aligned}
$$

Note that $\mathbb{E}(J_{\varepsilon_u}(\omega_j)J_{\varepsilon_v}(\omega_k)) = 0$, $1 \le j, k \le n/2$, and $\mathbb{E}(J_{\varepsilon_u}(\omega_j)\overline{J}_{\varepsilon_v}(\omega_k)) = 0$ if $|j-k| \ge p$ and

$$
\mathbb{E}(J_{\varepsilon_u}(\omega_j)\overline{J}_{\varepsilon_v}(\omega_k)) = \frac{\sigma_{uv}}{c_p}(-1)^{j-k}\binom{2p-2}{p-1+j-k}\mathbf{1}_{\{|j-k|<p\}}, \tag{53}
$$

where

$$c_p = \binom{2p-2}{p-1};$$

see [9]. Hence,

$$\begin{aligned}
&\mathbb{E}[(I_{\varepsilon\varepsilon,u_1v_1}(\omega_j) - \sigma_{u_1v_1})(I_{\varepsilon\varepsilon,u_2v_2}(\omega_k) - \sigma_{u_2v_2})] \\
&\quad = \mathbb{E}[I_{\varepsilon\varepsilon,u_1v_1}(\omega_j) I_{\varepsilon\varepsilon,u_2v_2}(\omega_k)] - \sigma_{u_1v_1}\sigma_{u_2v_2} \\
&\quad = \operatorname{cum}(J_{\varepsilon_{u_1}}(\omega_j), J_{\varepsilon_{u_2}}(\omega_k), \overline{J}_{\varepsilon_{v_1}}(\omega_j), \overline{J}_{\varepsilon_{v_2}}(\omega_k)) \\
&\qquad + \mathbb{E}(J_{\varepsilon_{u_1}}(\omega_j)\overline{J}_{\varepsilon_{v_2}}(\omega_k))\mathbb{E}(J_{\varepsilon_{u_2}}(\omega_j)\overline{J}_{\varepsilon_{v_1}}(\omega_k)) \\
&\quad = C\mathbf{1}_{\{|j-k|<p\}}
\end{aligned}$$

because (53) and the cumulant is 0 under Assumption 1. We have, by (5),

$$\begin{aligned}
&\mathbb{E}|S_{ab}(\omega_j)S_{ab}(\omega_k)| \\
&\quad \le C \sum_{u_1,u_2,v_1,v_2=1}^{q} (|\Psi_{au_1}(\omega_{\tilde{j}})||\overline{\Psi}_{au_2}(\omega_{\tilde{k}})||\overline{\Psi}_{bv_1}(\omega_{\tilde{j}})||\Psi_{bv_2}(\omega_{\tilde{k}})|)\mathbf{1}_{\{|j-k|<p\}} \\
&\quad \le C|1-e^{-i\omega_{\tilde{j}}}|^{-(d_{aa}+d_{bb})}|1-e^{-i\omega_{\tilde{k}}}|^{-(d_{aa}+d_{bb})}\mathbf{1}_{\{|j-k|<p\}}. \qquad \square
\end{aligned}$$

LEMMA 20. *Under Assumptions 1 and 2,*

$$b'\mathbf{A}\mathbf{f}(\omega_j)\mathbf{A}'b - b'\mathbf{A}_k\mathbf{f}_{kk}(\omega_j)\mathbf{A}_k'b = O_p(\omega_j^{-2d_k}(j^{d_k-d_{k-1}} + \omega_j^{d_k-d_{k+1}}))$$

*and*

$$b'\mathbf{A}_k\mathbf{f}_{kk}(\omega_j)\mathbf{A}_k'b - \mathcal{G}\omega_j^{-2d_k} = O_p(\omega_j^{-2d_k+\rho})$$

*for $1 \le j \le m_n$. Furthermore, there exists a constant $C$ such that*

$$b'\mathbf{A}\mathbf{f}(\omega_j)\mathbf{A}'b \ge C\omega_j^{-2d_k}(1-\varepsilon_k)$$

*for $1 \le j \le m_n$ and two constants, $C_1$ and $C_2$, such that*

$$C_1 > \mathcal{G} \ge C_2(1-\varepsilon_k),$$

*where $\varepsilon_k = O_p(n^{-\alpha_k})$.*

PROOF. Since, by Lemma 16,

$$b'\mathbf{A}_h\mathbf{f}_{h\ell}(\omega_j)\mathbf{A}_\ell'b = \begin{cases} O_p(\omega_j^{-2d_k}j^{2d_k-d_h-d_\ell}), & h<k, \ell \le k, \\ O_p(\omega_j^{-d_h-d_\ell}), & h,\ell > k, \\ O_p(\omega_j^{-d_k-d_\ell}j^{d_k-d_h}), & h \le k, \ell > k, \end{cases} \tag{54}$$

we have

$$\begin{aligned}
b'\mathbf{A}\mathbf{f}(\omega_j)\mathbf{A}'b &= b'\mathbf{A}_k\mathbf{f}_{kk}(\omega_j)\mathbf{A}_k'b + \sum_{h=0}^{s}\sum_{\substack{\ell=0\\ \ell\neq k}}^{s} b'\mathbf{A}_h\mathbf{f}_{h\ell}(\omega_j)\mathbf{A}_\ell' b \\
&= b'\mathbf{A}_k\mathbf{f}_{kk}(\omega_j)\mathbf{A}_k'b + O_p(\omega_j^{-2d_k} j^{d_k-d_{k-1}} + \omega_j^{-2d_{k+1}} + \omega_j^{-d_k-d_{k+1}}) \\
&= b'\mathbf{A}_k\mathbf{f}_{kk}(\omega_j)\mathbf{A}_k'b + O_p(\omega_j^{-2d_k}(j^{d_k-d_{k-1}} + \omega_j^{d_k-d_{k+1}})).
\end{aligned}$$

Since (7) and Assumption 2 imply that $\mathbf{f}_{kk}(\omega) = \mathbf{f}_{kk}^{\dagger}(0)\omega^{-2d_k} + O(\omega^{-2d_k+\rho})$ as $\omega \to 0$, we have, by Lemma 16,

$$\begin{aligned}
b'\mathbf{A}\mathbf{f}_{kk}(\omega_j)\mathbf{A}_k b' &= b'\mathbf{A}_k\mathbf{f}_{kk}^{\dagger}(0)\mathbf{A}_k' b\omega_j^{-2d_k} + O_p(\|b'\mathbf{A}_k\|^2\omega_j^{-2d_k+\rho}) \\
&= \mathcal{G}\omega_j^{-2d_k} + O_p(\omega_j^{-2d_k+\rho}).
\end{aligned}$$

We have shown the first two equations of the lemma. For the third equation, we have, by (54),

$$b'\mathbf{A}\mathbf{f}(\omega_j)\mathbf{A}'b = b'\left(\sum_{h,\ell=0}^{k}\mathbf{A}_h\mathbf{f}_{h\ell}(\omega_j)\mathbf{A}_\ell'\right)b + O_p(\omega_j^{-d_k-d_{k+1}}).$$

By Assumption 2 and Lemmas 16 and 17,

$$\begin{aligned}
b'\left(\sum_{h,\ell=0}^{k}\mathbf{A}_h\mathbf{f}_{h\ell}(\omega_j)\mathbf{A}_\ell'\right)b &= b'\left(\sum_{h,\ell=0}^{k}\omega_j^{-d_h-d_\ell}\mathbf{A}_h\mathbf{f}_{h\ell}^{\dagger}(0)\mathbf{A}_\ell'\right)b \\
&\quad + O_p\left(\omega_j^{-d_h+d_\ell+\rho}\sum_{h,\ell=0}^{k} b'\mathbf{A}_h\mathbf{A}_\ell' b\right) \\
&\geq \omega_j^{-2d_k}\lambda_{\min}\{\mathbf{f}^{\dagger}(0)\}\sum_{h,\ell=0}^{k} b'\mathbf{A}_h\mathbf{A}_\ell' b\mathbf{A}_h + O_p(\omega_j^{\rho}) \\
&\geq C\omega_j^{-2d_k}(1-\varepsilon_k) + O_p(\omega_j^{\rho}).
\end{aligned}$$

For the last inequality,

$$\mathcal{G} = b'\mathbf{A}_k\mathbf{f}_{kk}^{\dagger}(0)\mathbf{A}_k' b\omega_j^{-2d_k} \geq \omega_j^{-2d_k}\lambda_{\min}(\mathbf{f}_{kk}^{\dagger}(0))\|b'\mathbf{A}_k\|^2 \geq C\omega_j^{-2d_k}(1-\varepsilon_k),$$

by Lemma 17. The upper bound for $\mathcal{G}$ is due to the fact that

$$\mathcal{G} \leq \lambda_{\max}(\mathbf{f}_{kk}^{\dagger}(0))\|\mathbf{A}_k\|\|b\| = \lambda_{\max}(\mathbf{f}_{kk}^{\dagger}(0))\|\mathbf{A}_k\|. \qquad \square$$

LEMMA 21. *Let $\mathcal{L}_{m_n}(d)$, $\mathcal{M}_{m_n}(d)$ and $\mathcal{F}_{m_n}(d)$ be defined as in (21). Then if $d - d_k > -\frac{1}{2}$, there exists an $\varepsilon > 0$ such that*

(i) $$\mathcal{L}_{m_n}(d) = o_p(\omega_{m_n}^{2d-2d_k} m_n^{-\varepsilon}),$$

$$\text{(ii)} \qquad \mathcal{M}_{m_n}(d) = o_p(\omega_{m_n}^{2d-2d_k} m_n^{-\varepsilon}),$$

$$\text{(iii)} \qquad \mathcal{F}_{m_n}(d) - \mathcal{G}\omega_{m_n}^{2d-2d_k} = o_p(\omega_{m_n}^{2d-2d_k} m_n^{-\varepsilon}),$$

*under Assumptions 1 and 2.*

PROOF. We will only prove (ii); (i) and (iii) can be shown in a similar fashion using Lemmas 16, 18 and 20. Let $\mathbf{S}_{h\ell}(d)$ be the $(h,\ell)$th block matrix of $\mathbf{S}(d)$. By Lemma 19,

$$\begin{aligned}
&E\left\|\frac{1}{m_n}\sum_{j=1}^{m_n}\omega_{\tilde{j}}^{2d}\mathbf{S}_{h\ell}(\omega_j)\right\|^2 \\
&\quad = O\left(\frac{1}{m_n^2}\sum_{j=1}^{m_n}\sum_{k=j}^{p+j}\omega_{\tilde{j}}^{2d-d_h-d_\ell}\omega_{\tilde{k}}^{2d-d_h-d_\ell}\right) \\
&\quad = O\left(\frac{1}{m_n^2}\sum_{j=1}^{m_n}\sum_{k=j}^{p+j}\omega_{\tilde{j}}^{4d-2d_h-2d_\ell}\right) \\
&\quad = \begin{cases} O(n^{2d_h+2d_\ell-4d} m_n^{-2}\log m_n), & 4d-2d_h-2d_\ell \le -1, \\ O(\omega_{m_n}^{4d-2d_h-2d_\ell} m_n^{-1}), & 4d-2d_h-2d_\ell > -1. \end{cases}
\end{aligned}$$

Hence, we have

$$\begin{aligned}
&\left\|\frac{1}{m_n}\sum_{j=1}^{m_n}\omega_{\tilde{j}}^{2d}\mathbf{S}_{h\ell}(\omega_j)\right\| \\
&\quad = \begin{cases} O_p(n^{d_h+d_\ell-2d} m_n^{-1}\log^{1/2} m_n), & 2d-d_h-d_\ell \le -1/2, \\ O_p(\omega_{m_n}^{2d-d_h-d_\ell} m_n^{-1/2}), & 2d-d_h-d_\ell > -1/2. \end{cases}
\end{aligned} \tag{55}$$

Let

$$\mathcal{M}_{m_n}^{(h,\ell)}(d) = b'\mathbf{A}_h\left(\frac{1}{m_n}\sum_{j=1}^{m_n}\omega_{\tilde{j}}^{2d}\mathbf{S}_{h\ell}(\omega_j)\right)\mathbf{A}_\ell' b.$$

By Lemma 19 and (55), we have, for $h,\ell < k$,

$$\begin{aligned}
&\mathcal{M}_{m_n}^{(h,\ell)}(d) \\
&\quad = \begin{cases} O_p(\omega_{m_n}^{2d-2d_k} m_n^{-1-2d+2d_k}\log^{1/2} m_n), & 2d-d_h-d_\ell \le -1/2, \\ o_p(\omega_{m_n}^{2d-2d_k} m_n^{-1/2+\varepsilon}), & 2d-d_h-d_\ell > -1/2, \end{cases}
\end{aligned} \tag{56}$$

where $\varepsilon > 0$. By the same lemma and (55), we have, for $h, \ell \geq k$,

$$\mathcal{M}_{m_n}^{(h,\ell)}(d) = O_p(\omega_{m_n}^{2d-d_h-d_\ell} m_n^{-1/2}) = O_p(\omega_{m_n}^{2d-2d_k} \omega_{m_n}^{2d_k-d_h-d_\ell} m_n^{-1/2}) \tag{57}$$

and for $h < k$, $\ell \geq k$,

$$\mathcal{M}_{m_n}^{(h,\ell)}(d) = \begin{cases} o_p(\omega_{m_n}^{2d-2d_k} m_n^{-1-2d+2d_k} \log^{1/2} m_n), & 2d - d_h - d_\ell \leq -1/2, \\ o_p(\omega_{m_n}^{2d-2d_k} \omega_{m_n}^{d_k-d_\ell} m_n^{-1/2+d_k-d_h}), & 2d - d_h - d_\ell > -1/2. \end{cases} \tag{58}$$

Hence,

$$\mathcal{M}_{m_n}(d) = \sum_{h=0}^{s} \sum_{\ell=0}^{s} \mathcal{M}_{m_n}^{(h,\ell)}(d) = o_p(\omega_{m_n}^{2d-2d_k} m_n^{-\varepsilon}),$$

since $2d_k - d_h - d_\ell > 0$ in (57) and $-1 - 2d + 2d_k < 0$ in (56) and (58). □

LEMMA 22. *Under the assumptions of Theorem 2, if $d - d_k > -\frac{1}{4}$, then*

$$\mathcal{L}_{m_n}(d) = o_p(\omega_{m_n}^{2d-2d_k} m_n^{-1/2-\varepsilon})$$

*and*

$$\mathcal{M}_{m_n}^{(h,\ell)}(d) = \begin{cases} O_p(\omega_{m_n}^{2d-2d_k} m_n^{-1/2}), & h = \ell = k, \\ o_p(\omega_{m_n}^{2d-2d_k} m_n^{-1/2-\varepsilon}), & \text{otherwise}, \end{cases}$$

*Furthermore, under the assumptions of Theorem 2,*

$$\mathcal{F}_{m_n}(d_k) - \mathcal{G} = O_p(m_n^{d_k-d_{k-1}} + \omega_{m_n}^{d_k-d_{k+1}}),$$

*where the $O_p(m_n^{d_k-d_{k-1}})$ term is vacuous if $k = 0$ and the $O_p(\omega_{m_n}^{d_k-d_{k+1}})$ term is vacuous if $k = s$.*

PROOF. This lemma is a corollary of Lemma 21. □

**Acknowledgments.** The authors thank three referees for comments that substantially improved the paper and they especially thank the Associate Editor for a very detailed reading of the paper and for providing extremely useful suggestions.

Texas A&M University
Department of Statistics
3143 TAMU
College Station, Texas 77843-3143
USA
E-mail: wchen@stat.tamu.edu

New York University
44 West 4th Street
New York, New York 10012
USA
E-mail: churvich@stern.nyu.edu